\documentclass[pdftex,11pt,a4paper,parskip=half,headsepline,bibtotocnumbered]{scrartcl}

\usepackage[english]{babel}
\usepackage[onehalfspacing]{setspace}
\usepackage[automark]{scrlayer-scrpage}
\usepackage[utf8]{inputenc}
\usepackage{latexsym}
\usepackage[T1]{fontenc}
\usepackage{amsmath}
\usepackage{amssymb}
\usepackage{amsfonts}
\usepackage{amstext}
\usepackage{pdfsync}
\usepackage{amsthm}
\usepackage{extarrows}
\usepackage{graphicx}
\usepackage{xcolor}
\usepackage{longtable}
\usepackage[longtable]{multirow} 
\usepackage{booktabs} 
\usepackage{listings}
\usepackage{natbib}
\usepackage{url}
\usepackage{hyperref}

\definecolor{TUGreen}{rgb}{0.517,0.721,0.094}
 \definecolor{middlegray}{rgb}{0.5,0.5,0.5}
 \definecolor{lightgray}{rgb}{0.8,0.8,0.8}
 \definecolor{orange}{HTML}{e88909}
 \definecolor{yac}{rgb}{0.6,0.6,0.1} 
 \definecolor{bblue}{rgb}{0,0.4470,0.7410}
 \definecolor{rred}{rgb}{0.8500,0.3250,0.0980}
 \definecolor{turquoise}{HTML}{47C2D6}
 \definecolor{ggreen}{HTML}{12895a}

\author{Marléne Baumeister, Marc Ditzhaus, Markus Pauly}
\title{Quantile-based MANOVA: A new tool for inferring multivariate data in factorial designs}
\makeatletter

\clearscrheadfoot
\ihead{\headmark}
\automark{section}
\ohead{Baumeister, Ditzhaus and Pauly}
\newtheorem{satz}{Theorem}[]

\newtheorem{lemma}[satz]{Lemma}
\newtheorem{prop}[satz]{Proposition}
\newtheorem{ass}[satz]{Assumption}

\newcommand{\cale}{\mathcal{E}}
\newcommand{\br}{\mathbb{R}}
\newcommand{\bn}{\mathbb{N}}

\newcommand{\E}{E}

\newcommand{\cov}{\text{Cov}}

\newcommand{\tr}{\text{tr}}

\begin{document}

\thispagestyle{empty}
\maketitle

\section*{Abstract}
Multivariate analysis-of-variance (MANOVA) is a well established tool to examine multivariate endpoints. While classical approaches depend on restrictive assumptions like normality and homogeneity, there is a recent trend to more general and flexible procedures. 
In this paper, we proceed on this path, but do not follow the typical mean-focused perspective. 
Instead we consider general quantiles, in particular the median, for a more robust multivariate analysis. 
The resulting methodology is applicable for all kind of factorial designs and shown to be asymptotically valid. 
 Our theoretical results are complemented by an extensive simulation study for small and moderate sample sizes. An illustrative data analysis is also presented.

\textbf{Keywords:} Efron's Bootstrap; Multivariate Analysis of Variance; Factorial Designs; Quantile-Based Analysis; Nonparametric Inference; Heteroscedasticity

\section{Introduction}
In various fields, e.g. biology, ecology, medicine, or psychology, several outcome variables are of simultaneous interest leading to multivariate data. 
For example, an ecologist may study the aggression against predators and the relative reproductive success (fitness) of birds grouped by sex and colour morph \citep[cf.][]{boerner_aggression_2009}. 
Other examples are psychological tests or different medical quantities, e.g. heart rate, blood pressure, weight, or height of a patient.
As pointed out by \citet{warne_primer_2014}, multivariate analysis-of-variance (MANOVA) is \textit{ \glqq one of the most common multivariate statistical procedures in the social science literature\grqq}. 
However, classical MANOVA \citep{bartlett_note_1939, dempster_significance_1960, lawley_generalization_1938, pillai_new_1955, wilks_sample_1946} relies on restrictive assumptions as normality and homogeneity of covariances. 
But the \textit{ \glqq normality assumption becomes quasi impossible to justify when moving from univariate to multivariate observations\grqq} \citep{konietschke_parametric_2015} and, similarly, homogeneity is often implausible. 
To overcome these, several remedies have been suggested for tackling at least one of both issues. 
Thereby solutions have been developed for specific layouts, e.g. one- or two-way \citep{krishnamoorthy_parametric_2010, zhang_two-way_2011, bathke_testing_2018, zhang_linear_2022} as well as for general factorial designs \citep{konietschke_parametric_2015, friedrich_mats_2018}.
As common in statistical inference, all these proposals focus on the expectation (vector) and thus infer means or contrasts thereof. 
For heavy-tail distributions and in case of outliers the mean is not the appropriate statistical estimand \citep[cf.][]{maronna_robust_2006}.
Therefore, the present paper aims to introduce
\renewcommand{\labelenumi}{\roman{enumi})}
\begin{enumerate}
	\item a robust, quantile-based counterpart to mean-based MANOVA procedures
	\item without assuming a specific distribution class (such as normality or sphericity) 
	\item while allowing for potential heterogeneity
	\item in the framework of general factorial designs, including, e.g., higher-way layouts.
\end{enumerate}
To achieve these aims, we extend the recently proposed QANOVA \citep[quantile-based analy\-sis-of-\-variance, ][]{ditzhaus_qanova_2021} approach for univariate endpoints to multivariate settings.
The QANOVA procedure is a powerful alternative to the commonly established quantile regression \cite[cf.][]{koenker_quantile_2001, koenker_goodness_1999} and allows \textit{ \glqq the simple incorporation of interaction effects without a loss in power \grqq} \citep{ditzhaus_qanova_2021}.
Thus, it remains to answer the question how to extend QANOVA as there exists several possibilities to define multivariate quantiles, see, e.g. \citet{small_survey_1990}, \citet{serfling_quantile_2002} and \citet{becker_robustness_2013}. 
For example, the R-Package \texttt{MNM} \citep{nordhausen_multivariate_2011} allows for MANOVA analyses based on the spatial median by \cite{oja_descriptive_1983} and its affine equivariant version, the Hettmansperger–Randles median. \texttt{MNM} covers tests equivalent to Hotelling's $T^2$-Test for more than two samples and tests regarding randomized block designs \citep{oja_multivariate_2010}. 

Different to \texttt{MNM}, our method is based on the vector of marginal quantiles \citep{babu_joint_1988}. The marginal quantiles have the advantage, that they are computational more efficient and easier to interpret. In particular, it allows 
for compatible post hoc analyses on univariate components with the QANOVA. 
The proposed method will be established within a fully heterogeneous model and can be used for any quantile, not only for the median. For proofing correctness of our method we employ refined results on empirical quantile processes \citep{vaart_weak_2000,ditzhaus_qanova_2021}, and combine them with strategies and ideas from \citet{friedrich_mats_2018} for mean-based MANOVA.

The paper is structured as follows. 
The model, the estimators for the population quantiles are presented in Section \ref{set-up} together with a brief introduction to general factorial designs. 
Section \ref{sec:methods} presents the statistical methods. First (Section \ref{sec:tests}), 
the test statistics are constructed and their mathematical foundation is explained. Thereafter, covariance estimators based on kernel density estimators \citep{nadaraya_non-parametric_1965}, bootstrapping \citep{efron_bootstrap_1979, chung_exact_2013} and an interval-based strategy \citep{price_estimating_2001} are proposed (Section \ref{sec:covariance}). Finally, a (group-wise) bootstrap strategy is considered (Section \ref{sec:resampling}) to estimate the unknown limit distribution of the test statistics. The corresponding proofs are given in the Supplement. To investigate the method's small sample properties and compare it with existing methods, an extensive simulation study is carried out in Section \ref{simu}. 
An illustrative data analysis of Egyptian skulls complements our investigation (Section \ref{dataanalysis}). 
  The paper closes with a discussion and an outlook.

\section{Motivation and Set-Up}\label{set-up}
We consider a general, multivariate model based on mutually independent $d$-dimensional observation vectors for individuals from $k$ different (sub-)groups, e.g., representing different treatments or different epochs of antique objects as in Section~\ref{dataanalysis}. 
In detail, the $j$th observations vector in group $i$ is denoted by 
\begin{align*}
	\mathbf{X}_{ij}=\left(X_{ij1},\dots,X_{ijd}\right)' \sim \mathbb{P}_{i},  \quad i\in\{1,\ldots,k\},\;j\in\{1,\ldots,n_i\}.
\end{align*}
Here, $\mathbb{P}_i$ denotes the joint distribution with corresponding multivariate distribution function $F_i$.
Furthermore, let $F_{i\ell}$ be the marginal distribution function of $X_{i1\ell}$ $,i\in\{1,\ldots,k\},$ $\ell\in\{1,\ldots,d\}$ with existing density function $f_{i\ell}$. 
The joint distribution function of two entries $X_{ijm}$ and $X_{ij\ell}$ is denoted by $F_{im\ell},$ $i\in\{1,\ldots,k\},\,m,\ell\in\{1,\ldots,d\}$.
Throughout this paper, we like to infer 
the vector of marginal quantiles \citep[cf.][]{babu_joint_1988} $\mathbf{q}_i=\left(q_{i1},\dots,q_{id}\right)'$, where 
\begin{align*}
	q_{i\ell}=F_{i\ell}^{-1}(p)=\inf\left\{t\in\br|F_{i\ell}(t)\ge p\right\},\quad i\in\{1,\ldots,k\},\,m,\ell\in\{1,\ldots,d\},
\end{align*}
for a pre-specified quantile level $p\in(0,1)$, e.g. $p=0.5$ for medians. Within this setting, we want to develop testing procedures for the general null hypothesis
\begin{align}\label{eqn:null}
	\mathcal{H}_0:\mathbf{H}\mathbf{q}=\mathbf{0}_r, \quad \mathbf{q}=(\mathbf{q}_1',\ldots,\mathbf{q}_k')',
\end{align}
where $\mathbf{H}\in\br^{r\times dk}$ is a contrast matrix, i.e. $\mathbf{H}\mathbf{1}_{dk}=\mathbf{0}_r$, and $\mathbf{1}_r$ and $\mathbf{0}_r$ are the $r$-dimensional vectors of $0$'s and $1$'s, respectively. 
The concrete choice of $\mathbf{H}$ depends on the underlying research question and is similar to classical mean-based MANOVA.
For example, the one-way MANOVA hypothesis of no group effect is obtained by selecting $\mathbf{H}=\mathbf{P}_k\otimes\mathbf{I}_d$: 
\begin{align*}
	\mathcal{H}_0:\,\left(\mathbf{P}_k\otimes\mathbf{I}_d\right)\mathbf{q}=\mathbf{0}_{dk}\Leftrightarrow\mathcal{H}_0:\,\mathbf{q}_1=\dots=\mathbf{q_k}.
\end{align*}
Turning to a two-way layout with factors $A$ having $a$ levels and $B$ possessing $b$ levels, we split up the group index $i$ into $i=(i_1,i_2)$ for $i_1\in\{1,\dots,a\}$ and $i_2\in\{1,\dots,b\}$ resulting in $k=a\cdot b$ (sub-)groups. In a more lucid way, the multivariate quantile $\mathbf{q}_{i}$ can be decomposed into a general effect $\mathbf{q}^{\mu}$, main effects $\mathbf{q}_{i_1}^\alpha$, $\mathbf{q}_{i_2}^\beta$ and an interaction effect $\mathbf{q}_{i_1i_2}^{\alpha\beta}$ as
\[\mathbf{q}_{i}=\mathbf{q}_{(i_1,i_2)}=\mathbf{q}^\mu+\mathbf{q}_{i_1}^\alpha+\mathbf{q}_{i_2}^\beta+\mathbf{q}_{i_1i_2}^{\alpha\beta},\]
assuming the usual side conditions $\sum_{i_1=1}^a\mathbf{q}_{i_1}^\alpha=\sum_{i_2=1}^b\mathbf{q}_{i_2}^\beta=\sum_{i_1=1}^a\mathbf{q}_{i_1i_2}^{\alpha\beta}=\sum_{i_2=1}^b\mathbf{q}_{i_1i_2}^{\alpha\beta}=\mathbf{0}_d$ to ensure identifiability. Then, null hypotheses for main and interaction effects can be formulated as follows:
\begin{itemize}
	\item $\mathcal{H}_0(A):\left(\mathbf{P}_a\otimes\frac{1}{b}\mathbf{J}_b\otimes\mathbf{I}_d\right)\mathbf{q}=\mathbf{0}_{dk}\Leftrightarrow\mathbf{\bar q}_{1\cdot}=\dots=\mathbf{\bar q}_{a\cdot}\Leftrightarrow\mathbf{q}_{i_1}^\alpha=0\,\forall \,i_1\in\{1,\dots,a\}$,
	\item $\mathcal{H}_0(B):\left(\frac{1}{a}\mathbf{J}_a\otimes\mathbf{P}_b\otimes\mathbf{I}_d\right)\mathbf{q}=\mathbf{0}_{dk}\Leftrightarrow\mathbf{\bar q}_{\cdot 1}=\dots=\mathbf{\bar q}_{\cdot b}\Leftrightarrow\mathbf{q}_{i_2}^\beta=0\,\forall\,i_2\in\{1,\dots,b\}$,
	\item $\mathcal{H}_0(AB):\left(\mathbf{P}_a\otimes\mathbf{P}_b\otimes\mathbf{I}_d\right)\mathbf{q}=\mathbf{0}_{dk}\Leftrightarrow\mathbf{\bar q}_{\cdot\cdot}-{\mathbf{\bar q}}_{i_1\cdot}-{\mathbf{\bar q}}_{\cdot i_2}+{\mathbf{q}}_{i_1i_2}\equiv\mathbf{0}_d\,\forall\,i_1\in\{1,\dots,a\},$ $ i_2\in\{1,\dots,b\}\Leftrightarrow{\mathbf{q}}_{i_1i_2}^{\alpha\beta}=0\,\forall\,i_1\in\{1,\dots,a\},\,i_2\in\{1,\dots,b\}$.
\end{itemize}
Here, $\otimes$ denotes the Kronecker product of matrices, $\mathbf{P}_d=\mathbf{I}_d-\frac{1}{d}\mathbf{J}_d$ is the $d$-dimensional centring matrix, $\mathbf{I}_d$ is the $d$-dimensional identity matrix and $\mathbf{J}_d=\mathbf{1}_d'\mathbf{1}_d$ is the $d\times d$ matrix consisting of $1$'s only.
The $\mathbf{\bar{q}}_{i_1\cdot},\;\mathbf{\bar{q}}_{\cdot i_2}$ and $\mathbf{\bar q}_{\cdot\cdot}$ are the means over the dotted indices.
Higher-way layouts and also hierarchically designs with nested factors can be incorporated in a similar way, see e.g. \cite{pauly_asymptotic_2015,friedrich_mats_2018,friedrich_resampling-based_2019} for mean-based testing strategies.

We like to stress that different matrices $\mathbf{H}$ can describe the same null hypotheses but may affect the statistic's outcome \citep{sattler_testing_2022}. 
 In the sequel we therefore follow the common practice to reformulate \eqref{eqn:null} as $\mathbf{Tq}=\mathbf{0}_{dk}$ with the unique projection matrix $\mathbf{T}=\mathbf{H}'\left(\mathbf{HH}'\right)^+\mathbf{H}\in\br^{dk\times dk}$, where $\mathbf{A}^+$ denotes the Moore-Penrose inverse of the matrix $\mathbf{A}$. It is easy to check that $\mathbf{H}$ and $\mathbf{T}$ indeed lead to the same null hypotheses while the matrix $\mathbf{T}$ has the advantage of being unique, symmetric and idempotent \citep{pauly_asymptotic_2015, konietschke_parametric_2015, friedrich_wild_2017, bathke_testing_2018}.

To infer \eqref{eqn:null} based on real-data, the marginal quantiles are estimated via the empirical quantiles
	\begin{align*}
	\hat q_{i\ell}=\hat F_{i\ell}^{-1}(p)=\inf\{t\in\br|\hat F_{i\ell}(t)\ge p\}=X_{\lceil n_ip\rceil:n_i}^{(i\ell)},\,i\in\{1,\dots,k\},\,\ell\in\{1,\dots,d\},
	\end{align*}
where $X_{1:n_i}^{(i\ell)}\le\dots\le X_{n_i:n_i}^{(i\ell)}$ are the order statistics of the $\ell$-th component within group $i$ and $\hat F_{i\ell}(t)=\frac{1}{n_i}\sum_{j=1}^{n_i}\mathbf{1}_{\{X_{ij\ell}\le t\}}$ is the respective marginal empirical distribution function evaluated at time $t\in\mathbb{R}$. 
Together with the pre-chosen contrast matrix and respective covariance estimators (discussed in the next section), they are used in quadratic form-type test statistics to infer \eqref{eqn:null}.

\section{Statistical Methods}\label{sec:methods}

\subsection{Construction of Tests}\label{sec:tests}
As we want to develop an (at least) asymptotically valid method, we first recall the central limit theorem for marginal quantiles \citep{babu_joint_1988}. It relies on the following two standard regularity assumptions on the sample sizes and the distribution functions. Here and subsequently, all limits are meant as $n\to \infty$.
\begin{ass}\label{non-vanishing}
	The groups do not vanish, i.e. $(n_i/n)\to\kappa_i>0$. 
\end{ass}
\begin{ass}\label{ass}
	Let $F_{i\ell}$ be continuously differentiable at $q_{i\ell}$ with positive derivative $f_{i\ell}(q_{i\ell}) > 0$ for every $\ell\in\{1,\dots,d\}$ and $i\in	\{1,\dots,k\}$.
\end{ass}
	\begin{prop}[Theorem 2.1 of \citet{babu_joint_1988}]\label{prop1}
	Under Assumptions \ref{non-vanishing} and \ref{ass}  we have convergence in distribution
	\begin{align}
		\sqrt{n}\left(\hat q_{i\ell}-q_{i\ell}\right)_{\ell\in\{1,\dots,d\}}\stackrel{d}{\longrightarrow}\mathbf{Z}_i,\quad i\in\{1,\ldots,k\},
	\end{align}
	where $\mathbf{Z}_i$ is a zero-mean, multivariate normal distributed random vector with covariance matrix $\mathbf{\Sigma}^{(i)} = (\mathbf{\Sigma}_{\ell m}^{(i)})_{\ell,m=1,\dots,d}$ given by the entries 
	\begin{align}\label{eqn:deg_Sigma}
		\mathbf{\Sigma}_{\ell m}^{(i)}=\begin{cases}
			\cfrac{1}{\kappa_i}\cfrac{1}{f_{i\ell}^2(q_{i\ell})}\left[p-p^2\right],\;\ell=m,\\
			\cfrac{1}{\kappa_i}\cfrac{1}{f_{i\ell}(q_{i\ell})f_{im}(q_{im})}\left[F_{i\ell m}\left(F_{i\ell}^{-1}(p),F_{im}^{-1}(p)\right)-p^2\right],\;\ell\not=m.
		\end{cases}
	\end{align}
\end{prop}
For ease of convenience we present an empirical processes  based proof for Proposition~\ref{prop1} in the supplement. 
Here, a central limit theorem for multivariate quantiles is deduced from the functional delta-method for empirical processes \citep[cf.][Thm. 3.9.4]{vaart_weak_2000}.  
In principle, Proposition~\ref{prop1} and the group's independence allow us to construct quadratic form test statistics %
 in terms of the vector $\sqrt{n}(\mathbf{\widehat q} - \mathbf{q})$. 
For this purpose, we only require an appropriate estimator for the unknown limiting covariance matrix $\mathbf{\Sigma}:=\bigoplus_{i=1}^k\mathbf{\Sigma}^{(i)}$.
Consistent proposals are discussed in Section~\ref{sec:covariance}. Let us suppose for a moment, that $\mathbf{\widehat\Sigma}$ consistently estimates $\mathbf{\Sigma}$. 
Then we propose a so-called ANOVA-type statistic (ATS) \citep{brunner_box-type_1997} and a modified ANOVA-type statistic (MATS) \citep{friedrich_mats_2018} to infer $\mathbf{Tq}=\mathbf{0}_{dk}$:
\begin{align}
	\mathbf{ATS}_n(\mathbf{T})=n\frac{(\mathbf{T\hat q)'\mathbf{T\hat q}}}{\tr\left(\mathbf{T\hat\Sigma T}\right)}, \quad \mathbf{MATS}_n(\mathbf{T})=n(\mathbf{T\hat q)'\left(\mathbf{T}\hat\Sigma_0\mathbf{T}\right)^+\mathbf{T\hat q}}.
\end{align}
Here, $\mathbf{\Sigma}_0$ and $\mathbf{\hat\Sigma}_0$ denote the matrices containing only the diagonal elements of $\mathbf{\Sigma}$ and $\mathbf{\hat\Sigma}$, respectively. 
As described in  \citet{sattler_testing_2022}, both test statistics can be unified into the following general form 
\begin{align}\label{teststatistic}
	S_n(\mathbf{T})=n(\mathbf{T\hat q})'E\left(\mathbf{T},\mathbf{\hat\Sigma}\right)\mathbf{T\hat q},
\end{align}
where $E(\mathbf{T},\mathbf{\hat\Sigma})=\tr(\mathbf{T\hat\Sigma T})^{-1}\mathbf{I}_{dk}$ for the ATS and $E(\mathbf{T},\mathbf{\hat\Sigma})=(\mathbf{T}\mathbf{\hat\Sigma}_0\mathbf{T})^+$ for the MATS.
In the supplement, we prove that both versions are consistent for $E(\mathbf{T},\mathbf{\Sigma})$.
The following theorem summarizes the asymptotic distribution of $S_n(\mathbf{T})$.
	\begin{satz}\label{distributiontests}
 Let $\mathbf{\hat\Sigma}$ be a consistent covariance matrix estimator of $\mathbf{\Sigma}$ and assume that Assumptions \ref{non-vanishing} and \ref{ass} hold.
	\begin{enumerate}
		\item Under $\mathcal{H}_0:\,\mathbf{Tq}=\mathbf{0}_{dk}$, the test statistic $S_n$ converges in distribution to a weighted sum of $\chi_1^2$ distributed random variables, i.e.
		\begin{align}\label{formdistrtests}
			S_n(\mathbf{T})\stackrel{d}{\longrightarrow}B=\sum_{i=1}^{dk}\lambda_iB_i,
		\end{align}
		where $B_i\stackrel{iid}{\sim}\chi_1^2$ and $\lambda_i\ge 0,\,i\in\{1,\dots,dk\},$ are the eigenvalues of $E\left(\mathbf{T},\mathbf{\Sigma}\right)\mathbf{T\Sigma T}$.
		\item Under $\mathcal H_1:\,\mathbf{Tq}\neq\mathbf{0}_{dk}$, $S_n$ converges in probability to $\infty$.
	\end{enumerate}
	\end{satz}
Let $b_\alpha$ be the $(1-\alpha)$-quantile of $B$ in \eqref{formdistrtests}. 
From Theorem \ref{distributiontests} we can deduce that the test  $\varphi_n=\mathbf{1}{\{S_n(\mathbf{T})>b_\alpha\}}$ is of asymptotic level $\alpha$ for $\mathcal{H}_0:\,\mathbf{Tq}=\mathbf{0}$. 
Furthermore, it is consistent for any alternative $\mathcal H_1:\,\mathbf{Tq}\neq\mathbf{0}_{dk}$.  
However, the distribution of $B$ depends on unknown parameters through $\mathbf{\Sigma}$. 
Thus, $b_\alpha$ is in general unknown and we consider a bootstrap procedure to approximate it, which we discuss in Section~\ref{sec:resampling}. 
But first we first we address the pending question regarding to the estimation of $\mathbf{\Sigma}$.

\subsection{Estimation of the Covariance Matrix}\label{sec:covariance}
Variance estimation of the median or general quantiles is not easy in the univariate setting and various strategies can be found in the literature \citep{maritz_note_1978, mckean_comparison_1984, bonett_confidence_2006, chung_exact_2013}. 
At a first glance, the situation becomes even more complicated in the multivariate set-up. 
However, a careful observation of \eqref{eqn:deg_Sigma} yields a simple relationship between the diagonal and off-diagonal elements:
\begin{align*}
	\mathbf{\Sigma}^{(i)}_{\ell m}=\sqrt{\mathbf{\Sigma}^{(i)}_{\ell\ell}\mathbf{\Sigma}^{(i)}_{mm}}\,\frac{F_{i\ell m}(q_{i\ell},q_{im})-p^2}{p-p^2},\,\ell\not=m.
\end{align*}
Thus, the known univariate strategies to estimate the variances $\mathbf{\Sigma}^{(i)}_{\ell\ell}$, $\ell\in\{1,\ldots, k\}$, can be combined with an estimator for $F_{i\ell m}(q_{i\ell},q_{im})$. 
For the latter, let us first introduce the joint empirical distribution function $\hat F_{i\ell m}$ defined by
\begin{align*}
\hat F_{i\ell m}(t_1,t_2)=\frac{1}{n_i}\sum_{j=1}^{n_i}\mathbf{1}_{(-\infty,t_1]\times(-\infty,t_2]}(X_{ij\ell},X_{ijm})
\end{align*}
for $\;i\in\{1,\dots,k\}$,  $\ell,m\in\{1,\dots,d\}$ and $t_1,t_2\in\mathbb{R}$.
Under the following regularity assumption, consistency of $\hat F_{iab}(\hat q_{ia},\hat q_{ib})$ for $F_{iab}(q_{ia},q_{ib})$ holds:
	\begin{ass}\label{jointcontinuous}
	For every $i\in\{1,\dots,k\}$ and every $\ell,m\in\{1,\dots,d\}$, the joint distribution function $F_{i\ell m}$ is continuous at $(q_{i\ell},q_{im})=(F_{i\ell}^{-1}(p),F_{im}^{-1}	(p))$. 
	\end{ass}
	\begin{prop}[Theorem 2.2 of \cite{babu_joint_1988}]\label{jointconsistent}
	Under Assumption \ref{jointcontinuous} and \ref{non-vanishing},  the estimator $\hat F_{i\ell m}(\hat q_{i\ell},\hat q_{im})$ converges in probability to $F_{i\ell m}(q_{i\ell},q_{im})$. 
	\end{prop}
Consequently, we obtain a general form of estimators for $\mathbf{\Sigma}^{(i)}$:
\begin{align*}
	\mathbf{\hat\Sigma}_{\ell m}^{(i)}=\begin{cases}
		\cfrac{n}{n_i}\,\hat\sigma^2_{i\ell}(p),\,\ell=m,\\
		\cfrac{n}{n_i}\,\hat\sigma_{i\ell}(p)\hat\sigma_{im}(p)\,\cfrac{\hat F_{i\ell m}(\hat q_{i\ell},\hat q_{im})-p^2}{p-p^2},\,\ell\not=m,
	\end{cases}
\end{align*}
where $\hat\sigma_{i\ell}^2(p)$ is a consistent estimator for the asymptotic variance $\sigma_{i\ell}^2(p)=\kappa_i \mathbf{\Sigma}_{\ell\ell}^{(i)}$ 
of the marginal, centred empirical quantiles $\sqrt{n}(\widehat q_{i\ell} - q_{i\ell})$. 
We follow \citet{ditzhaus_qanova_2021} and study three different choices for $\hat\sigma_{i\ell}^2(p)$ considered in the univariate case. All approaches produce consistent estimators under the respective assumptions for $\mathbf{\hat\Sigma}^{(i)}$ \citep[cf.][]{ditzhaus_qanova_2021}.

\subsubsection{Kernel Estimator}
The \textit{kernel-based approach} uses the strong consistent kernel density estimator by \citet{nadaraya_non-parametric_1965} to estimate the densities $f_{i\ell}$. It is given by 
	\begin{align}\label{kernel density}
	\hat f_{K,i,\ell}(x)=\frac{1}{n_ih_{n_{i}\ell}}\sum_{j=1}^{n_i}K_{i\ell}	\left(\frac{x-X_{ij\ell}}{h_{n_i\ell}}\right),\,i\in\{1,\dots,k\};\,\ell\in\{1,\dots,d\},
	\end{align}
where $K_{i\ell}$ is a kernel and $h_{n_i\ell}$ is a bandwidth, $i\in\{1,\dots,k\}$; $\ell\in\{1,\dots,d\}$.
For its strong consistency we require: 
\begin{ass}\label{kernelassumption} 
		Suppose for every $i\in\{1,\dots,k\};\ell\in\{1,\dots,d\}$ that $K_{i\ell}$ is of bounded variation, $f_i$ is uniformly continuous and the series $\sum_{m=1}^\infty\exp(-\gamma mh_{n_i\ell}^2)$ converges for every choice of $\gamma>0$.
	\end{ass}
This leads to the following consistent estimator for $\sigma_{i\ell}^2(p)$:
	\begin{align*}
	\hat\sigma^{2}_{i\ell,K}(p)=\frac{1}{\hat f_{K,i,\ell}^2(\hat q_{i\ell})}(p-p^2).
	\end{align*}
	
\subsubsection{Bootstrap Estimator}\label{sec:bootstrap}
The \textit{bootstrap approach} was originally proposed by \cite{chung_exact_2013}, who borrowed the idea from \cite{efron_bootstrap_1979}. To introduce it, consider the bootstrap samples $X_{i1\ell}^*,\dots,X_{in_i\ell}^*$, $i\in\{1,\ldots,k\}$, $\ell\in\{1,\ldots,d\},$ drawn mutually independent and with replacement from the observations $X_{i1\ell},\dots,X_{in_i\ell}$. 
We denote all estimators based on the bootstrap sample by a $^*$, e.g. $\hat q^*$. 
Then, the bootstrap sample quantile estimator $\hat q_{i\ell}^*$ can be calculated and its conditional mean squared error, given data, 
is given by 
\begin{align}\label{bootstrap}
		\left(\hat\sigma_{i\ell}^*(p)\right)^2=n_i\sum_{j=1}^{n_i}\left(X_{j:n_i}^{(i\ell)}-\hat q_{i\ell}\right)^2\underbrace{P^*\left(X_{\lceil n_ip\rceil:n_i}^{(i\ell)*}=X_{j:n_i}^{(i\ell)}|\mathbf{X}_{i\ell}\right)}_{:=P^*_{ij\ell}}.
		\end{align}
As explained by \citet{efron_bootstrap_1979}, $P^*_{ij\ell}$ can be rewritten for every $\ell\in\{1,\dots,d\}$ as
		\begin{align*}
		P^*_{ij\ell}=P\left(B_{n_i,\frac{j-1}{n_i}}\le\lceil n_ip\rceil-1\right)-P\left(B_{n_i,\frac{j}{n_i}}\le\lceil n_ip\rceil-1\right),
		\end{align*}
	where $B_{n,p}$ denotes a binomial distributed random variable with size parameter $n$ and success probability $p$.
\citet{ghosh_note_1984} proved that the estimator $\left(\hat\sigma_{i\ell}^*(p)\right)^2$ is consistent for $\sigma_{i\ell}^2(p)$ under the following moment condition.
		\begin{ass}\label{bootstrapassumption}
For some $\delta>0$ we have $\max_{i\in\{1,\dots,k\},\,\ell\in\{1,\dots,d\}}\E\left(|X_{i1\ell}|^\delta\right)<\infty$.
		\end{ass}
		
\subsubsection{Interval-based Estimator}	
An \textit{interval-based approach} was initially suggested by \citet{mckean_comparison_1984} and later modified by \citet{price_estimating_2001} for the median.
The methodology can easily be adapted to handle general quantiles \citep{bonett_confidence_2006,ditzhaus_qanova_2021}.
 The principle idea is to start with the asymptotic confidence interval $(X_{l_i(p):n_i}^{(i\ell)},X_{l_i(p):n_i}^{(i\ell)})$ for $q_{i\ell}$. Its length (asymptotically) depends on the (unknown) standard deviation $\sigma_{i\ell}(p)$.  Basic calculation yields the following estimator:
		\begin{align}\label{interval-based}
		\hat\sigma_{i\ell,PB}^2(p)=\left(\sqrt{n_i}\;\cfrac{X_{u_i(p):n_i}^{(i\ell)}-X_{l_i(p):n_i}^{(i\ell)}}{2z_{\alpha_{n_i\ell}^*(p)/2}+2n_i^{-1/2}}\right)^2,
		\end{align}
	where $l_i(p)=\max\{1,\lfloor n_ip-z_{\alpha/2}\sqrt{n_ip(1-p)}\rfloor\}$ is the lower and $u_i(p)=\min\{n_i,$ $\lfloor n_ip+$ $z_{\alpha/2}$ $\sqrt{n_ip(1-p)}\rfloor\}$  is the upper limit of the binomial interval, and $z_{\alpha/2}$ denotes the $\left(1-{\alpha}/{2}\right)$-quantile of the standard normal distribution.
Note that $l_i(p)$ and $u_i(p)$ are independent of the dimension $\ell\in\{1,\dots,d\}$, and typically $\alpha=0.05$ is chosen for their computation.
{\citet{price_estimating_2001} set
		\begin{align*}
		\alpha_{n_i\ell}^*(p)=
		\begin{cases}
		\alpha_{n_i\ell}(p)=1-\displaystyle\sum_{j=l_i(p)+1}^{u_i(p)-1}\begin{pmatrix}n_i\\j\end{pmatrix}p^j(1-p)^{n_i-j},\;n_i\le100,\\
		z_{0.05/2}\approx1.96,\;n_i>100
		\end{cases}
		\end{align*}
in the denominator of $\hat\sigma_{i\ell,PB}^{2}(p)$. 
The case distinction is motivated by the computation time of $\alpha_{n_i\ell}(p)$ which becomes quite demanding for larger sample sizes. Since we have $\alpha_{n_i\ell}(p)\longrightarrow\alpha$ by the central limit theorem, the use of $\alpha_{n_i\ell}(p)$ is only necessary for small to moderate sample sizes.}

\subsection{Bootstrapping the Test Statistic}\label{sec:resampling}

Having one of the presented estimators for the covariance matrices at hand, we are able to calculate the respective ATS or MATS.
However, the limiting distribution of $B$ from Theorem \ref{distributiontests} remains unknown and the $(1-\alpha)$-quantile $b_\alpha$ of $B$ cannot be computed.
That is why we consider a group-wise nonparametric bootstrap approach to approximate it.
This resampling strategy was already applied by \citet{friedrich_mats_2018} for their mean-based MANOVA procedure and is known to be effective for various testing problems \citep{konietschke_parametric_2015, dobler_nonparametric_2020, liu2020resampling}. 
As in Section~\ref{sec:bootstrap},  we consider a $d$-dimensional bootstrap sample $\{\mathbf{X}_{i1}^*,\dots,\mathbf{X}_{in_i}^*\}$ drawn with replacement from the original observation vectors $\{\mathbf{X}_{i1},\dots,\mathbf{X}_{in_i}\}$ for every group $i\in\{1,\dots,k\}$. 
Moreover, we add a $^*$ to all statistics which are calculated from the bootstrap sample, e.g. $\mathbf{\hat q^*}$ denotes the bootstrap quantile vector. 
Observe that under the null hypothesis $\mathcal H_0:\mathbf{T}\mathbf{q}=\mathbf{0}_{dk}$ the test statistic can be written as 
\begin{align*}
	S_n(\mathbf{T})=n\left[\mathbf{T}\left(\mathbf{\hat q}-\mathbf{ q}\right)\right]'E\left(\mathbf{T},\mathbf{\hat\Sigma}\right)\mathbf{T}\left(\mathbf{\hat q}-\mathbf{q}\right).
\end{align*}
For its bootstrap counterpart, the estimators $\mathbf{\hat q}$ and $\mathbf{\hat\Sigma}$ are replaced by their bootstrap versions $\mathbf{\hat q^*}$ and $\mathbf{\hat\Sigma^*}$  and the unknown quantile vector $\mathbf{q}$ is substituted by its the empirical counterpart $\mathbf{\hat q}$. 
Consequently, we obtain the bootstrap statistic
\begin{align}\label{bootstrapteststat}
	S_n^*(\mathbf{T})=n\left[\mathbf{T}\left(\mathbf{\hat q}^*-\mathbf{\hat q}\right)\right]'E\left(\mathbf{T},\mathbf{\hat\Sigma}^*\right)\mathbf{T}\left(\mathbf{\hat q}^*-\mathbf{\hat q}\right).
\end{align}
{Ich w{\"u}rde stattdessen schreiben: We can derive the (conditional) asymptotic behaviour of $S_n^*(T)$ by slightly adopting the argumentation for Proposition~\ref{prop1} and combine that with the bootstrap results of \citet{vaart_weak_2000}} {to get an equivalent result to Theorem \ref{distributiontests}:} 
	\begin{satz}\label{bootdistributiontest}
	Let $\mathbf{\hat \Sigma}$ be a consistent estimator for $\mathbf{\Sigma}$ and let $\mathbf{\hat\Sigma}^*$ denote its consistent bootstrap version.
	Then, the bootstrap test statistic $S_n^*(\mathbf{T})$ given in (\ref{bootstrapteststat}) converges always, conditionally given the data, in distribution to a real-valued random variable $B^*$, i.e. we have under $\mathcal{H}_0:\,\mathbf{Tq}=\mathbf{0}_{dk}$ as well as under $\mathcal{H}_1:\,\mathbf{Tq}\not=\mathbf{0}_{dk}$ 
		\begin{align*}
		\sup_{x\in\br}\left|P\left(S_n^*(\mathbf{T})\le x|\,\mathbf{X}\right)-P\left(B^*\le x\right)\right|\stackrel{P}{\longrightarrow}0.
		\end{align*}
	Hereby, the distribution of $B^*$ depends on the underlying setting and can be expressed by $B^*=\sum_{i=1}^{dk}\lambda^*_iB_i$ where $\lambda^*_i\geq 0$ and $B_i\sim\chi_1^2$.
	Under $\mathcal{H}_0$, the distribution of $B^*$ coincides with the limit null distribution of $S_n(\mathbf{T})$, i.e. $B\stackrel{d}{=}B^*$.
	\end{satz}	
Our proposed resampling test use $b_\alpha^*(\mathbf{X})$, the empirical $(1-\alpha)$-quantile of the conditional distribution function $x\mapsto P(S_n^*(\mathbf{T})\le x|\,\mathbf{X})$ as critical value.
This leads to the test $\varphi_n^*=\mathbf{1}{\{S_n(\mathbf{T})>b_\alpha^*(\mathbf{X})\}}$.
Under $\mathcal{H}_0$, Theorem~\ref{bootdistributiontest} implies that $b_\alpha^*(\mathbf{X})$ converges in probability to $b_\alpha$ given the data.
Thus, combining Lemma 1 of \citet{janssen_how_2003}, Theorem \ref{distributiontests} and Theorem~\ref{bootdistributiontest}, we obtain that the resampling test $\varphi_n^*$ is asymptotically exact, i.e. $\E_{\mathcal{H}_0}(\varphi_n^*)\to\alpha$. 
Moreover, $\varphi_n^*$ is even consistent for general alternatives $\mathcal{H}_0:\mathbf{Tq}\not=\mathbf{0}_{dk}$. 
To accept this, we first deduce from Theorem \ref{bootdistributiontest} that $b_\alpha^*(\mathbf{X})$ converges in probability to some $\widetilde b\in\mathbb{R}$ given the data under $\mathcal H_1$. Combining this observation with Theorem \ref{distributiontests} ii) and Theorem 7 of \citet{janssen_how_2003} yields the desired consistency.
	
\section{Simulations}\label{simu}
\begin{figure}[ht!]
\includegraphics[scale=0.62]{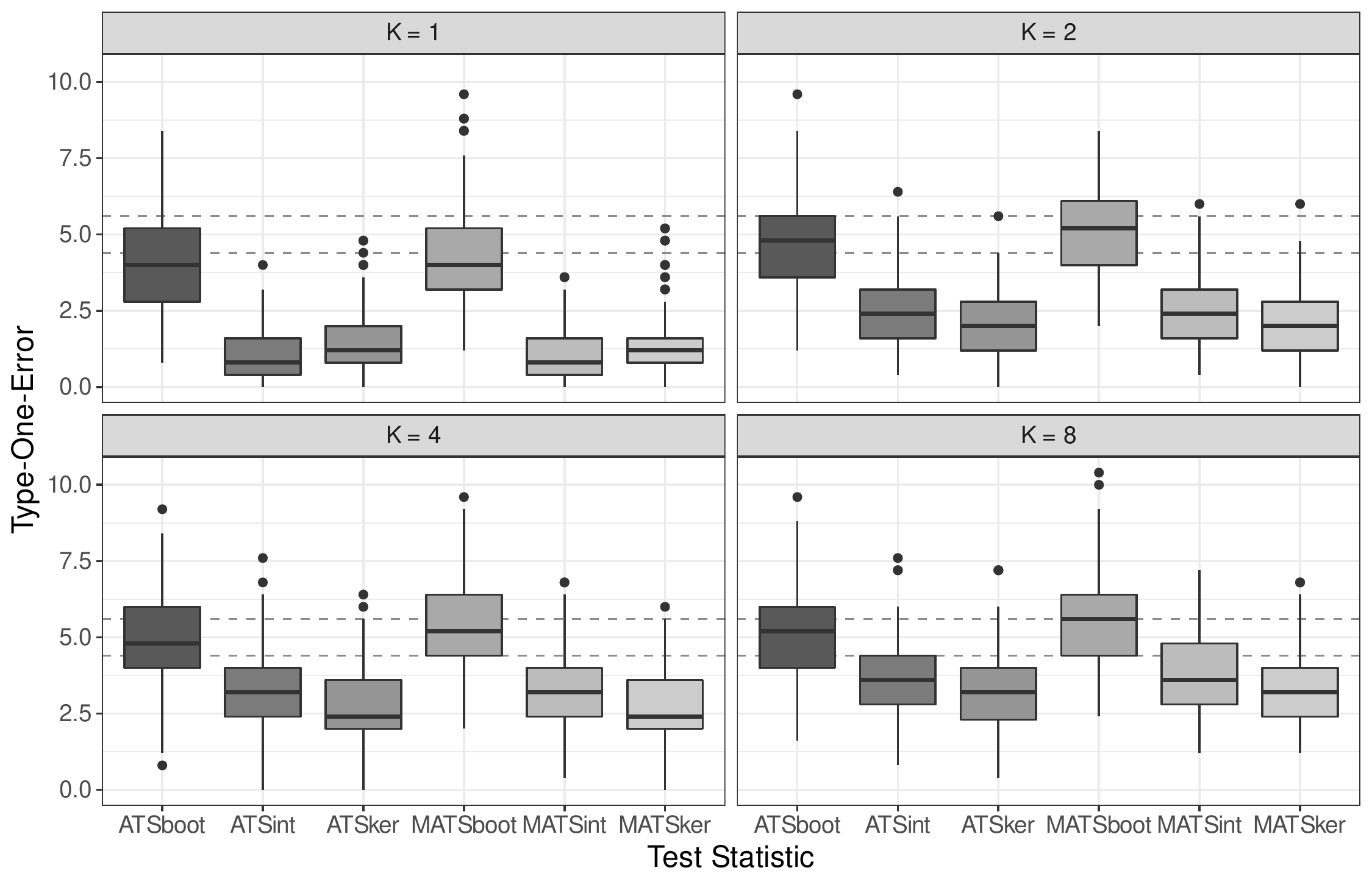}
\caption{{Type I error in $\%$ of the six versions of the QMANOVA divided into two types (ATS, MATS) and three covariance estimator (boot, int, ker). 
The results are divided by the sample size factor $K\in\{1,2,4,8\}$, i.e. $\mathbf{n}=K\mathbf{n}^{(r)}$, $r=1,2,3$.
}}
\label{graph:boxplots}
\end{figure}
To assess the tests' performances for small and moderate sample sizes, we conducted a simulation study. 
The idea of the data generation is as follows. 
We generate median-centred data $\mathbf{e}_{ij}=\left(e_{ij1},\dots,e_{ijn_i}\right)', \,e_{ij\ell}=Y_{ij\ell}-\text{median}(Y_{ij\ell})$ and choose for $Y_{ij\ell}$ the following five different distributions to cover symmetric and skewed scenarios: (a) the standard normal distribution $Y_{ij\ell}^{(1)}\sim N_{0,1}$, (b) the Student's $t$-distribution with $2$ degrees of freedom $Y_{ij\ell}^{(2)}\sim t_2$, (c) the Student's t-distribution with $3$ degrees of freedom $Y_{ij\ell}^{(3)}\sim t_3$, (d) the standard log-normal distribution $Y_{ij\ell}^{(4)}\sim LN_{0,1}$ and (e) the Chi-square distribution with $3$ degrees of freedom $Y_{ij\ell}^{(5)}\sim \chi^2_3$.
Certain homoscedastic and hetereoscedastic covariance settings are realized by multiplying the square root of different covariance matrices to this data.
Furthermore, we considered six different covariance matrices representing homoscedastic and heteroscedastic scenarios, which are displayed below for $d=4$:	
{
\renewcommand{\labelenumi}{\arabic{enumi})}
	\begin{enumerate}
	\item $\mathbf{\Sigma}^{(1)}=\mathbf{I}_d-\frac{1}{2}\left(\mathbf{J}_d-\mathbf{I}_d\right)=\mathbf{\Sigma}^{(2)}$,
	\item $\mathbf{\Sigma}^{(1)}=\left(0.6^{|a-b|}\right)_{a,b=1}^d=\mathbf{\Sigma}^{(2)}$,
	\item $\mathbf{\Sigma}^{(1)}=\mathbf{I}_d+\frac{1}{2}\left(\mathbf{J}_d-\mathbf{I}_d\right),\,\mathbf{\Sigma}^{(2)}=3\,\mathbf{I}_d+\frac{1}{2}\left(\mathbf{J}_d-\mathbf{I}_d\right)$,
	\item $\mathbf{\Sigma}^{(1)}=\left(0.6^{|a-b|}\right)_{a,b=1}^d,\,\mathbf{\Sigma}^{(2)}=\left(0.6^{|a-b|}\right)_{a,b=1}^d+2\,\mathbf{I}_d$,
	\item $\mathbf{\Sigma}^{(1)}=\begin{pmatrix}
	1 & 0.6 & 0.36 & 0.18\\
	0.6 & 1 & 0.6 & 0.3\\
	0.36 & 0.6 & 1 & 0.5\\
	0.18 & 0.3 & 0.5 & 0.25
	\end{pmatrix},\,\mathbf{\Sigma}^{(2)}=\mathbf{\Sigma}^{(1)}+0.5\,\mathbf{J}_d$,
	\item $\mathbf{\Sigma}^{(1)}=\begin{pmatrix}
	1 & 0 & 0 & 0\\
	0 & \sqrt{2} & 0 & 0\\
	0 & 0 & 2 & 1\\
	0 & 0 & 1 & 0.5
	\end{pmatrix},\,\mathbf{\Sigma}^{(2)}=\mathbf{\Sigma}^{(1)}+0.5\,\mathbf{J}_d$.
	\end{enumerate}}
{Here, the fifth covariance setting is a modification form the second, where the elements $\sigma_{1dd}$, $\sigma_{1(d-1)d}$ and $\sigma_{1d(d-1)}$ of $\mathbf{\Sigma}^{(1)}$  are modified as described.
The sixth setting is based on $\text{diag}(\sqrt{2^s})$ for $s\in\{0,\dots,d-1\}$, where the $d$-th row and column is replaced by half the row or rather the column before.
To create data which is median centred, we calculate the empirical median $\mathbf{M}_{ij}$ of $(\mathbf{\Sigma}^{(i)})^{\frac{1}{2}}\mathbf{e}_{ij}$ from an extra sample with the size $n=10^7$ and withdraw $\mathbf{M}_{ij}$ from the data.
Therefore, our simulated data can be described by the following model:}
{
\begin{align*}
\mathbf{X}_{ij}=(\mathbf{\Sigma}^{(i)})^{\frac{1}{2}}\mathbf{e}_{ij}-\mathbf{M}_{ij}\sim F_i,\,i\in\{1,\dots,k\},j\in\{1,\dots,n_i\}.
\end{align*}}	
The aforementioned data generating process and the choice of the different settings is adapted from the one-way layout simulation in \citet[Sec.~5]{friedrich_mats_2018}.
To consider small and large sample scenarios, we chose balanced and unbalanced small samples $\mathbf{n}^{(1)}=(10,10)',\,\mathbf{n}^{(2)}=(10,20)'$ and $\mathbf{n}^{(3)}=(20,10)'$ as well as its multiples $K\cdot\mathbf{n}^{(r)}$ for $K\in \{2,4,8\}$.
As a benchmark, we compare our method with the mean-based resampling MATS proposed by \citet{friedrich_mats_2018}.
This method is implemented in the R-package \texttt{MANOVA.RM} 
 \citep{friedrich_manovarm_2022} as the functions \texttt{MANOVA()} or \texttt{MANOVA.wide()} (same tests for different data formats).
We simulated two versions of the mean-based MATS, one is characterized by a parametric bootstrap and the other by a wild bootstrap with Rademacher weights, which are both implemented in the aforementioned package. 
All simulations are calculated with the computing environment R \citep{r_core_team_r_2021}, Version 4.0.0, for \texttt{nsim = 5000} simulation runs and \texttt{nboot = 2000} bootstrap iterations.
As  in \citet{ditzhaus_qanova_2021}, we used the classical Gaussian kernel for the kernel density estimation and calculate the bandwidth $h$ by Silverman's rule-of-thumb \citep[Eq. 3.31]{silverman_density_1998} with the R-function \texttt{bw.rnd0()}.

\subsection{Type I Error}
\begin{figure}[h!]
\includegraphics[scale=0.58]{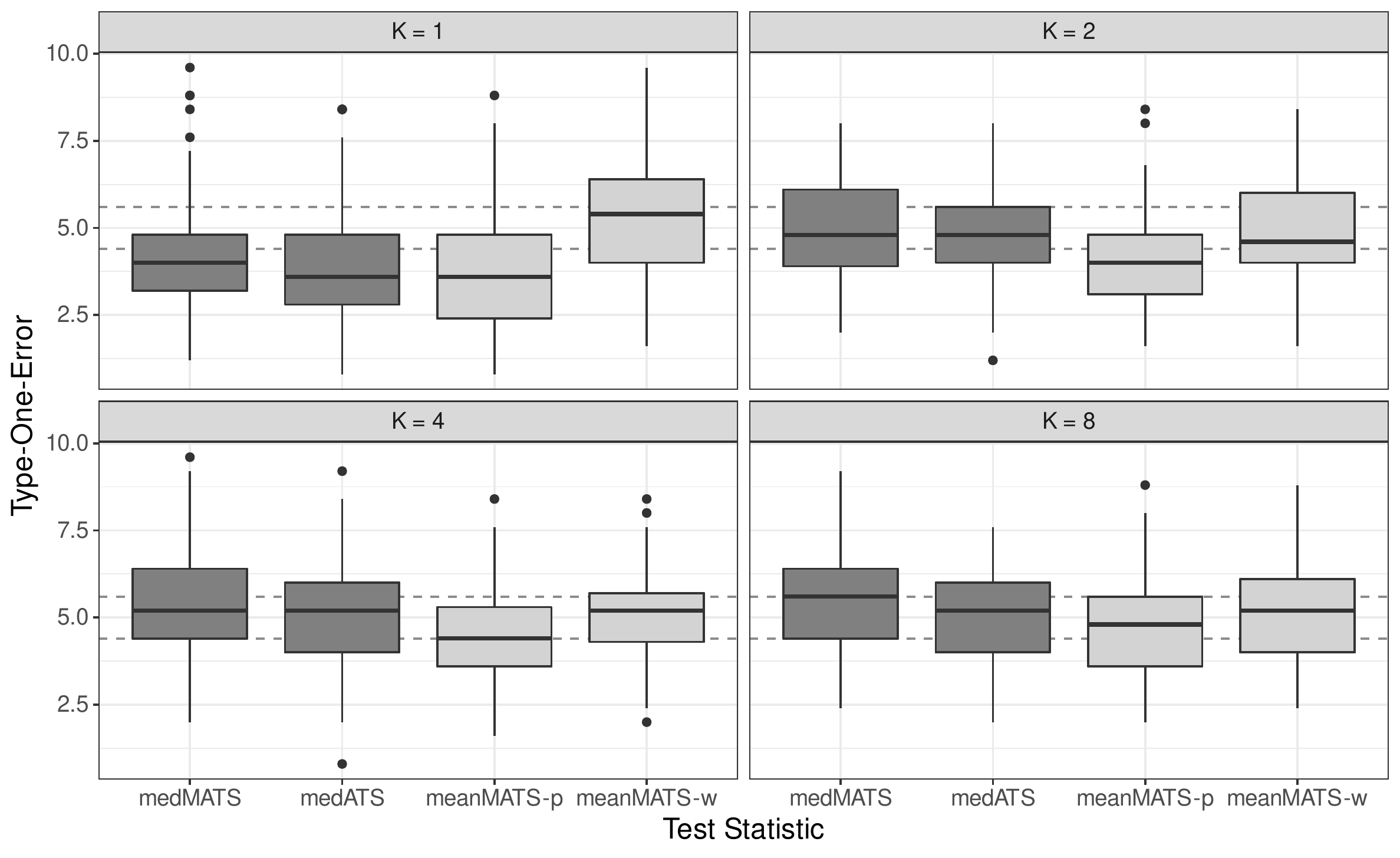}
\caption{Type I error in $\%$ of the QMANOVA with bootstrap covariance estimator and the mean-based MATS by \citet{friedrich_mats_2018} with parametric (p) and wild (w) bootstrap for the one-way layout, the symmetric distributions ($N_{0,1}$, $t_2$, $t_3$), all covariance settings 1)- 6) as well as for balanced and unbalanced sample sizes $\mathbf{n}= K\cdot\mathbf{n}^{(r)}$ for $K\in\{1,2,4,8\}$.}
\label{graph:comp}
\end{figure}
In this subsection, we discuss the type I error control of all procedures in a one-way layout and present further results for a $2\times 2$-design in the supplement. 
In detail, we considered a multivariate set-up with $k=2$ groups and $d=\{4,8\}$ dimensions. Moreover, we restricted to the median $\mathbf{m}_i,\,i\in\{1,2\}$, ($p=0.5$), because it is the most relevant quantile for statistical analysis and is comparable to the mean. 
This lead us to the null hypothesis  $\mathcal{H}_0:\mathbf{m}_1=\mathbf{m}_2$ for the layout matrix $\mathbf{T}=\mathbf{P}_2\otimes\mathbf{I}_4$ and in all to $720$ different scenarios.
In Figure \ref{graph:boxplots}, the different tests are named by a combination of the used test statistic and its covariance estimator.
It is apparent that the bootstrap covariance estimator has the best performance regarding the type I error control.
Overall, the MATS and the ATS test statistic have a similar type-one-error control. However, the MATS with the bootstrap covariance estimator performs the best as one can observe from the close position of the boxes to the binomial interval $[4.4,5.6]$.
With the other covariance estimators, ATS and MATS show a quite conservative type I error control.
In general, the observed type I error rates comes closer to the 5\%-benchmark for larger $K$. 
This is in line with the theoretical findings from Section \ref{sec:tests}. 
All in all, we can only recommend the ATS and MATS combined with the bootstrap covariance estimator.

In the next step, we compare the two favorable QMANOVA methods, from now denoted by \textit{medMATS} and \textit{medATS}, with the mean-based MATS of \cite{friedrich_mats_2018} denoted by \textit{meanMATS-p} and \textit{meanMATS-w} to differentiate between the parametric (p) and wild (w) bootstrap versions.
For a fair and appropriate comparison, we restrict ourselves to the symmetric distributions such that the mean-based hypothesis $\mathcal{H}_0:\boldsymbol{\mu}_1=\boldsymbol{\mu}_2$ and the median-based hypothesis $\mathcal{H}_0:\mathbf{m}_1=\mathbf{m}_2$ coincide.
The type I error rates for all four tests are summarized in Table~\ref{table:tenten} for $\mathbf{n}^{(1)}=(10,10)$ and in the supplement for $\mathbf{n}^{(2)}$ and $\mathbf{n}^{(3)}$.
For completeness reasons, we also include the error rates for two QMANOVA strategies for the nonsymmetric distributions $(LN_{0,1},\chi^2_3$). 
The results inside the binomial interval $[4.4,5.6]$ for the significance level $\alpha=5\%$ are printed in bold. 
{A detailed study of Table \ref{table:tenten} exhibits that the MATS performs slightly better than the ATS since the simulated type I errors are more frequent in $[4.4,5.6]$ for the MATS (17 times) than for the ATS (14 times).}
To further judge the performance for larger sample sizes, we summarized the type I error rates for all symmetric distributions and all different sample sizes $K\mathbf{n}^{(r)}$, $r=1,2,3$, in boxplots displayed in Figure \ref{graph:comp} divided into the different choices of the scaling factor $K\in\{1,2,4,8\}$.
In the Figure \ref{graph:comp}, the medMATS tends to be more conservative (empirical type I error smaller than $4.4$) and the meanMATS-w tends to be more liberal (empirical type I error larger than $5.6$).
In case of small sample sizes $(K=1)$, the medATS and both MATS approaches exhibit a conservative type I error control for almost all settings. 
{In settings with $d=8$ and the covariance choices 4) or 6), a conservative type I error control can even be found for larger samples $(K\in\{2,4,8\})$, see Table~\ref{table:tenten}.}
Additionally, a liberal behaviour most often occurs for the covariance setting 5) in combination with the dimension $d=4$ (Table~\ref{table:tenten}). 
For further details on the influence of the simulation scenarios' aspects we refer to the Supplement. There, we explain that the choice of the covariance setting  influences the performance most wile the QMANOVA method is mostly robust.
Moreover, we also present type I error simulation study for a two-way layout in the supplement. 
There we simulated scenarios for a $2\times 2$ design with  dimension $d=4$ for similar covariance settings. The corresponding findings are similar to the one-layout and thus omitted here.

\begin{table}
\centering
\small
\singlespacing
\begin{tabular}{llllllllll}
\toprule
\multirow[c]{3}{*}{$\Sigma$} & \multirow[c]{3}{*}{Distr} & \multicolumn{4}{c}{$d=4$}  & \multicolumn{4}{c}{$d=8$}\\ \cmidrule(l{0.5em}r{0.5em}){3-6} \cmidrule(l{0.5em}r{0.5em}){7-10}
& &  \multicolumn{2}{c}{median}  & \multicolumn{2}{c}{meanMATS}& \multicolumn{2}{c}{median}  & \multicolumn{2}{c}{meanMATS}\\ \cmidrule(l{0.5em}r{0.5em}){3-4} \cmidrule(l{0.5em}r{0.5em}){5-6} \cmidrule(l{0.5em}r{0.5em}){7-8} \cmidrule(l{0.5em}r{0.5em}){9-10}
 & &  MATS & ATS & param & wild & MATS & ATS & param & wild\\
  \midrule
\multirow[t]{6}{*}{1} & $N_{0,1}$  & 2.4 & 2.0 & 3.2 & 4.0 & 2.4 & 2.8 & 3.2 & 4.0 \\ 
  & $t_2$  & 1.6 & $\mathbf{4.4}$ & 2.0 & 4.0 & 4.0 & 6.0 & 4.0 & 6.8 \\ 
  & $t_3$ & $\mathbf{4.8}$ & $\mathbf{5.6}$ & 4.0 & $\mathbf{4.4}$ & $\mathbf{4.4}$ & 3.6 & 3.6 & $\mathbf{4.4}$ \\ 
  & $LN_{0,1}$ & 3.2 & 3.6 & - & - & $\mathbf{4.8}$ & $\mathbf{4.4}$ & - & - \\ 
 & $\chi^2_3$ & 2.0 & 2.0 & - & - & $\mathbf{4.8}$ & 4.0 & - & - \\ 
  \addlinespace
 \multirow[t]{6}{*}{2} & $N_{0,1}$ & $\mathbf{4.8}$ & 3.2 & 3.2 & $\mathbf{4.4}$ & $\mathbf{4.4}$ & 2.8 & 2.8 & 3.2 \\ 
  & $t_2$ & $\mathbf{4.4}$ & $\mathbf{4.4}$ & $\mathbf{4.8}$ & 7.6 & $\mathbf{4.8}$ & $\mathbf{4.8}$ & 2.8 & $\mathbf{5.2}$ \\ 
 & $t_3$ & $\mathbf{4.4}$ & 3.6 & 3.6 & $\mathbf{4.4}$ & 1.2 & $\mathbf{4.4}$ & 3.6 & $\mathbf{4.4}$ \\ 
 & $LN_{0,1}$ & 6.0 & 6.4 & - & - & 3.6 & $\mathbf{5.2}$ & - & - \\ 
  & $\chi^2_3$ & $\mathbf{5.2}$ & $\mathbf{4.4}$ & - & - & 3.6 & 4.0 & - & - \\ 
  \addlinespace
  \multirow[t]{6}{*}{3} & $N_{0,1}$  & 3.6 & 1.2 & 3.2 & 4.0 & 3.2 & 2.8 & $\mathbf{4.8}$ & 7.2 \\ 
 & $t_2$  & 4.0 & $\mathbf{4.8}$ & 3.2 & 6.4 & 3.6 & 3.2 & 3.2 & $\mathbf{5.2}$ \\ 
 & $t_3$  & 2.8 & 2.8 & 3.6 & $\mathbf{5.6}$ & 3.2 & 2.4 & 3.6 & 6.8 \\ 
  & $LN_{0,1}$  & 1.6 & 2.0 & - & - & 3.2 & 4.0 & - & - \\ 
  & $\chi^2_3$ & 6.0 & 3.6 & - & - & 4.0 & 2.4 & - & - \\ 
  \addlinespace
\multirow[t]{6}{*}{4} & $N_{0,1}$ & 1.6 & 1.6 & 2.4 & 3.6 & 2.8 & 2.8 & $\mathbf{4.8}$ & 6.4 \\ 
 & $t_2$ & 4.0 & 2.4 & 2.4 & $\mathbf{5.6}$ & 2.8 & 3.6 & 1.2 & 3.6 \\ 
  & $t_3$  & 3.6 & 3.2 & 1.2 & 1.6 & 2.8 & 2.0 & 2.8 & $\mathbf{5.6}$ \\ 
  & $LN_{0,1}$ & $\mathbf{4.8}$ & $\mathbf{4.8}$ & - & - & 3.6 & 2.0 & - & - \\ 
 &  $\chi^2_3$ & $\mathbf{5.2}$ & $\mathbf{5.6}$ & - & - & 2.4 & 1.6 & - & - \\ 
  \addlinespace
 \multirow[t]{6}{*}{5} & $N_{0,1}$ & 8.8 & 8.4 & 8.8 & 9.6 & $\mathbf{4.4}$ & $\mathbf{4.8}$ & $\mathbf{4.8}$ & 6.0 \\ 
 & $t_2$ & 3.6 & 3.6 & $\mathbf{5.6}$ & 6.0 & 1.2 & 2.4 & 2.8 & 3.6 \\ 
 & $t_3$  & 7.6 & 8.4 & 6.0 & 6.0 & 4.0 & 3.2 & $\mathbf{5.2}$ & 6.0 \\ 
 & $LN_{0,1}$ & 3.2 & 3.6 & - & - & 6.4 & 8.0 & - & - \\ 
 &  $\chi^2_3$ & 8.4 & 6.4 & - & - & $\mathbf{4.8}$ & 4.0 & - & - \\ 
  \addlinespace
 \multirow[t]{6}{*}{6} & $N_{0,1}$  & 7.2 & 6.4 & $\mathbf{4.8}$ & $\mathbf{5.6}$ & 1.6 & 2.4 & 1.6 & 2.8 \\ 
 & $t_2$ & $\mathbf{4.8}$ & 3.6 & 1.6 & 4.0 & 2.4 & 3.2 & 2.0 & 4.0 \\ 
 & $t_3$ & 4.0 & $\mathbf{4.4}$ & 3.2 & $\mathbf{4.4}$ & 2.8 & $\mathbf{4.4}$ & 2.4 & 6.0 \\ 
 & $LN_{0,1}$ & $\mathbf{5.6}$ & 6.8 & - & - & 3.2 & 2.8 & - & - \\ 
  & $\chi^2_3$ & $\mathbf{5.2}$ & 3.2 & - & - & 4.0 & 6.4 & - & - \\ 
   \bottomrule
   \end{tabular}
   \caption{Type-I error rate in $\%$ (nominal level $\alpha = 5\%$) for testing $\mathcal{H}_0: \mathbf{m}_1=\mathbf{m}_2$ in an one-way layout. We present all settings with the sample size $(10,\,10)$ of MATS and ATS combined with the bootstrap covariance estimator.
The two comparing meanMATS models are named with \textit{meanMATS-p} for a parametric bootstrap and with \textit{meanMATS-w} for the wild bootstrap approach.
The results inside the binomial interval for $\alpha$ $[4.4,5.6]$ are printed in bold.}
\label{table:tenten}
\end{table}

\subsection{Power}
\begin{figure}[ht!]
\includegraphics[scale=0.643]{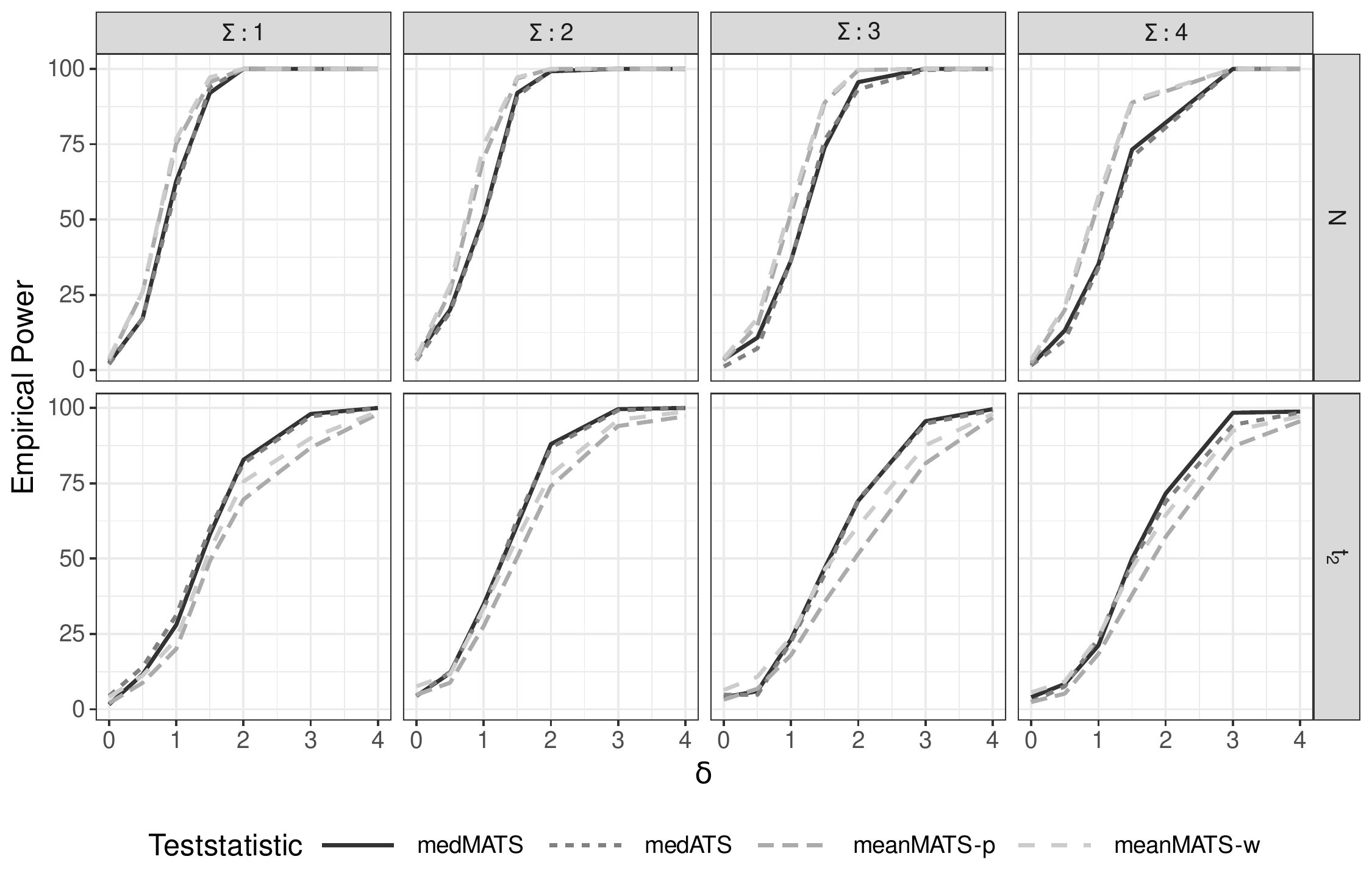}
\caption{Empirical Power in Percent of the QMANOVA and the mean-based MATS by \citet{friedrich_mats_2018} for normal distributed (N) and $t_2$-distributed ($t_2$) data, the small sample size $(10,10)$ and the covariance settings $1-4$.}
\label{graph:power}
\end{figure}
For the power comparison, we restrict to a representative subset of the settings from the previous section, namely to $d=4$, the covariance settings 1) - 4), the samples $Kn^{(1)}$ and $Kn^{(3)}$ for $K=2,4$ while we still consider all five distributions. 
To obtain alternative settings, we shift the data from the first group by $\delta\in\{0.5,1,1.5,2,3,4\}$, i.e. $\mathbf{X}_{1j}=\delta+\left((\mathbf{\Sigma}^{(1)})^{\frac{1}{2}}\mathbf{e}_{1j}-\mathbf{M}_{1j}\right),j\in\{1,\dots,n_i\}$
Moreover, we restrict to a representative subset of the settings, namely to $d=4$, the covariance settings 1) - 4), the samples $Kn^{(1)}$ and $Kn^{(3)}$ for $K=2,4$ while we still consider all five distributions.
We again compare the two favorable QMANOVA methods with the mean-based MATS for symmetric distributions.
Figure \ref{graph:power} includes the empirical power of the four methods for the normal and the $t_2$-distribution and all simulated covariance settings.
Studying Figure \ref{graph:power} one can observe, that the mean-based tests are more powerful in the cases with the normal distribution, but for the $t_2$-distribution we can see the exact opposite.
This observation fits to the power simulation results in \citet{ditzhaus_qanova_2021}.
There, an explanation for this is also given: mean and median are as location estimators asymptotically different efficient in the distributional scenarios.
The sample median is the better location estimator in case of heavy-tailed data like the $t_2$-distribution, but for normal distributed data the situation is reversed.
For all distributions the power increases faster for bigger sample sizes and again, the unbalanced designs do not seem to have any influence on this. 

\section{Illustrative Data Analysis}\label{dataanalysis}

\begin{figure}
\includegraphics[scale=0.64]{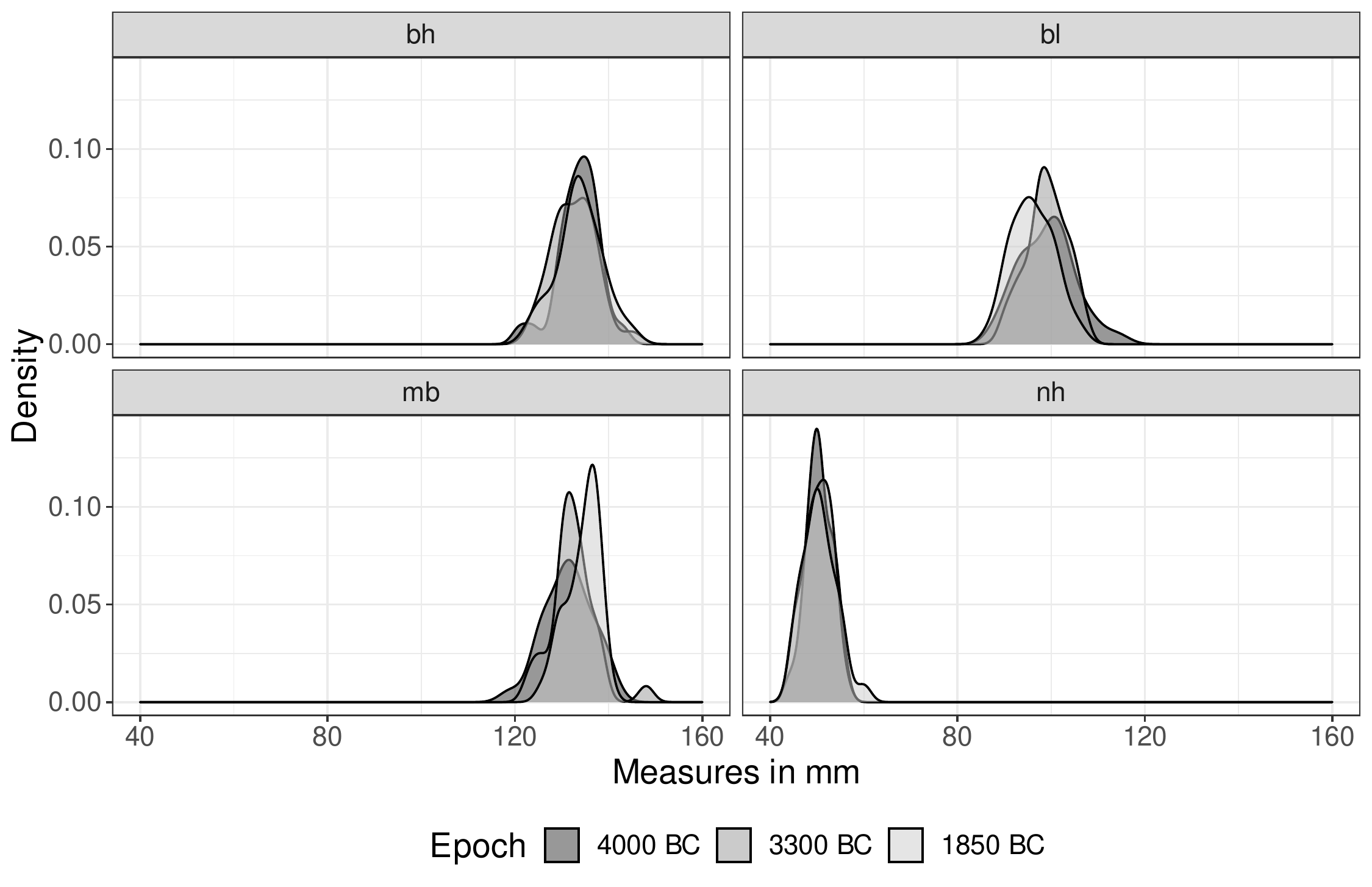}
\caption{Kernel density estimates \citep{nadaraya_non-parametric_1965} based on Gaussian kernels of the four ($d=4$) characteristics of the skulls mb (maximal breadth), bh (basibregmatic height), bl (basialveolar length) and nh (nasal height) divided by three epochs ($k=3$) arround 4000 BC, around 3300 BC and around 1850 BC. The bandwidth of the kernel densities is chosen by Silverman's rule-of-thumb \citep[Eq. 3.31]{silverman_density_1998}.}
\label{graph:skulls}
\end{figure}

To illustrate the new methods on real data, we re-analyse the Egyptian skulls data set from \citet{everitt_handbook_2006} available in the R-Package \texttt{HSAUR} \citep{everitt_hsaur_2022}.
For $90$ skulls, there are four variables ($d=4$) measured in mm and denoted by mb (maximal breadth), bh (basibregmatic height), bl (basialveolar length) and nh (nasal height).
The skulls can be divided into three groups ($k=3$) which are characterised by time periods in years around 4000 BC ($i=1$), around 3300 BC ($i=2$) and around 1850 BC ($i=3$) \citep{oja_multivariate_2010}.
The data is balanced with 30 skulls per group \citep{everitt_hsaur_2022}.
All four measurements together characterize the skulls in their basic shape. 
We are interested in inferring whether there are differences between the three epochs. 
Figure~\ref{graph:skulls} shows kernel density plots for each univariate measurement. We observe that the data is rather heavy-tailed and potentially heteroscedastic in each of the four measures. 
Thus, a median-based approach is reasonable leading to the multivariate null hypothesis $\mathcal{H}_0^{123}:\,\mathbf{m}_{1}=\mathbf{m}_{2}=\mathbf{m}_{3}$.
Due to its convincing type I error control in our simulation study, we choose the quantile-based MATS combined with the bootstrap covariance estimator and compare it with the mean-based MATS \citep{friedrich_mats_2018} using wild bootstrap critical values.
Similar to the simulation study, all tests are computed based upon $2,000$ bootstrap iterations. 
The resulting $p$-values are given in Table~\ref{table:testresult}.
\begin{table}[t]
\centering
\begin{tabular}{lcccc}
\toprule
\multirow[c]{3}{*}{Test} & \multicolumn{4}{c}{Hypothesis}\\ 
\cmidrule{2-5}
& $\mathcal{H}_0^{123}:$ & $\mathcal{H}_{0}^{12}:$ & $\mathcal{H}_{0}^{23}:$ & $\mathcal{H}_{0}^{13}:$\\
& $\mathbf{m}_{1}=\mathbf{m}_{2}=\mathbf{m}_{3}$ & $\mathbf{m}_{1}=\mathbf{m}_{2}$ & $\mathbf{m}_{2}=\mathbf{m}_{3}$ & $\mathbf{m}_{1}=\mathbf{m}_{3}$\\
\midrule
medMATS & $0.038^*$ & $0.726$ & $0.042^*$ & $0.007^*$ \\
\addlinespace
meanMATS-w & $0.062$ & $0.756$ & $0.061$ & $0.038$\\
\bottomrule
\end{tabular}
\caption{$p$-values for the different null hypotheses and the tests \textit{medMATS} (quantile-based MATS with bootstrap covariance estimator) and \textit{meanMATS-w} (mean-based MATS with wild bootstrap resampling). 
The $p$-values marked with a ${}^*$ indicate the rejected hypootheses at multiple level $5\%$.
}
\label{table:testresult}

\end{table}
It can be seen, that the \textit{medMATS} rejects the null hypothesis at significance level $0.05$ whereas the mean-based \textit{meanMATS-w} does not. 
A probable reason for this is that the data is nearly symmetric and heavy-tailed in all dimensions (cf. Figure \ref{graph:skulls}). 
In such settings the \textit{medMATS} exhibit a better power performance compared to \textit{meanMATS-w}.
After rejecting this global null hypothesis, it is intuitive to perform group-wise post hoc analyses.
Consequently, one can formulate all pairs hypotheses: $\mathcal{H}_{0}^{12}: \mathbf{m}_{1}=\mathbf{m}_{2}$, $\mathcal{H}_{0}^{23}: \mathbf{m}_{2}=\mathbf{m}_{3}$ and $\mathcal{H}_{0}^{13}: \mathbf{m}_{1}=\mathbf{m}_{3}$.  
The tests' p-values are displayed in Table \ref{table:testresult}.
Altogether, the considered hypotheses form a closed testing procedure \citep{gabriel_simultaneous_1969} and thus do not need adjustment. That is why the multiple \textit{meanMATS-w} and \textit{medMATS} tests control the family-wise-error-rate without an extra adjustment of the $p$-values. We obtain different test results from both methods:
The multiple \textit{medMATS} test rejects the global null and detects a difference in the 3300 BC and the 1850 BC skulls at the 5\% level. 
In comparison, the multiple mean-based \textit{meanMATS-w} test does not reject the global null hypothesis $\mathcal{H}_{0}^{123}$ at the 5\% level. Thus, it would also not detect the difference between the 3300 BC and the 1850 BC skulls as no pairwise posthoc comparisons would have been performed. 


\section{Conclusion and Outlook}\label{concl}
We have introduced six statistical tests in a general MANOVA-set-up (QMANOVA) regarding marginal quantiles . 
These are based on different test statistics: an ANOVA-type statistic (ATS) and a modified ATS (MATS), both in combination with three different covariance estimators (based on a kernel, bootstrap and interval-based approach). All statistics are quadratic forms in the normalized vector of pooled quantiles and can be seen as a generalization of the univariate QANOVA presented in \citet{ditzhaus_qanova_2021}. 
As all test statistics are no asymptotic pivots, we propose a non-parametric bootstrap approach to calculate critical values. We analyze the corresponding limit behavior and prove that the resulting tests are asymptotic exact and consistent.
In fact, the methods are asymptotically valid 
in general factorial designs and do not postulate homoscedasticity or a specific distribution. 
In an extensive simulation study focusing on the median, it turned out that the MATS with the bootstrap covariance estimator performs the best among the six proposed QMANOVA methods. In particular, it is robust against various aspects of data and performs well on heavy-tailed and skewed data with equal, unequal and singular covariance structures. 
We additionally compared its performance with the mean-based MATS proposed by \citet{friedrich_mats_2018} as benchmark method. 
In-line with theoretical properties of means and medians, our power simulation study showed that the median-based QMANOVA performs better than the corresponding mean-based approach in terms of power. 
This has also been confirmed in an illustrative data analyses with
symmetric and heavy-tailed data in a three-way MANOVA setting. The test results suggest that using the QMANOVA instead of an mean-based method is an added value in this application.


Apart from the illustrative data analyses, the focus of the paper was on deriving global test procedures. Having rejected a global null, post-hoc analyses on the components or factor levels are of interest and multiplicity may become an issue. In the three-way MANOVA setting of the data example this was no issue. But it would be for more complex situtations. Thus, we plan to derive multiple contrast tests (MCTs) for contrasts of marginal medians and quantiles in the future. Here, concepts from \citet{gunawardana_nonparametric_2019} could be adapted. As the derived methods can directly be inverted in confidence regions we would also like to derive simultaneous confidence regions and intervals that are compatible to the MCTs decisions. This would allow a deeper insight into the behaviour of the estimates. Moreover, similar to mean-based MANCOVA \citep{zimmermann2020multivariate}, we plan to derive QMANCOVAs that allow for covariate adjustments.  



\section{Acknowledgement}
This work has been partly supported by the Research Center Trustworthy Data Science and Security (\url{https://rc-trust.ai}), one of the Research Alliance centers within the \href{https://uaruhr.de}{UA Ruhr}.

\bibliographystyle{apalike}
\bibliography{Literatur, MarkusMarc}

\appendix
\newpage
\section{Additional Simulation Results}
\subsection{Additional Results regarding the One-Way Layout}
\begin{figure}[h]
\includegraphics[scale=0.57]{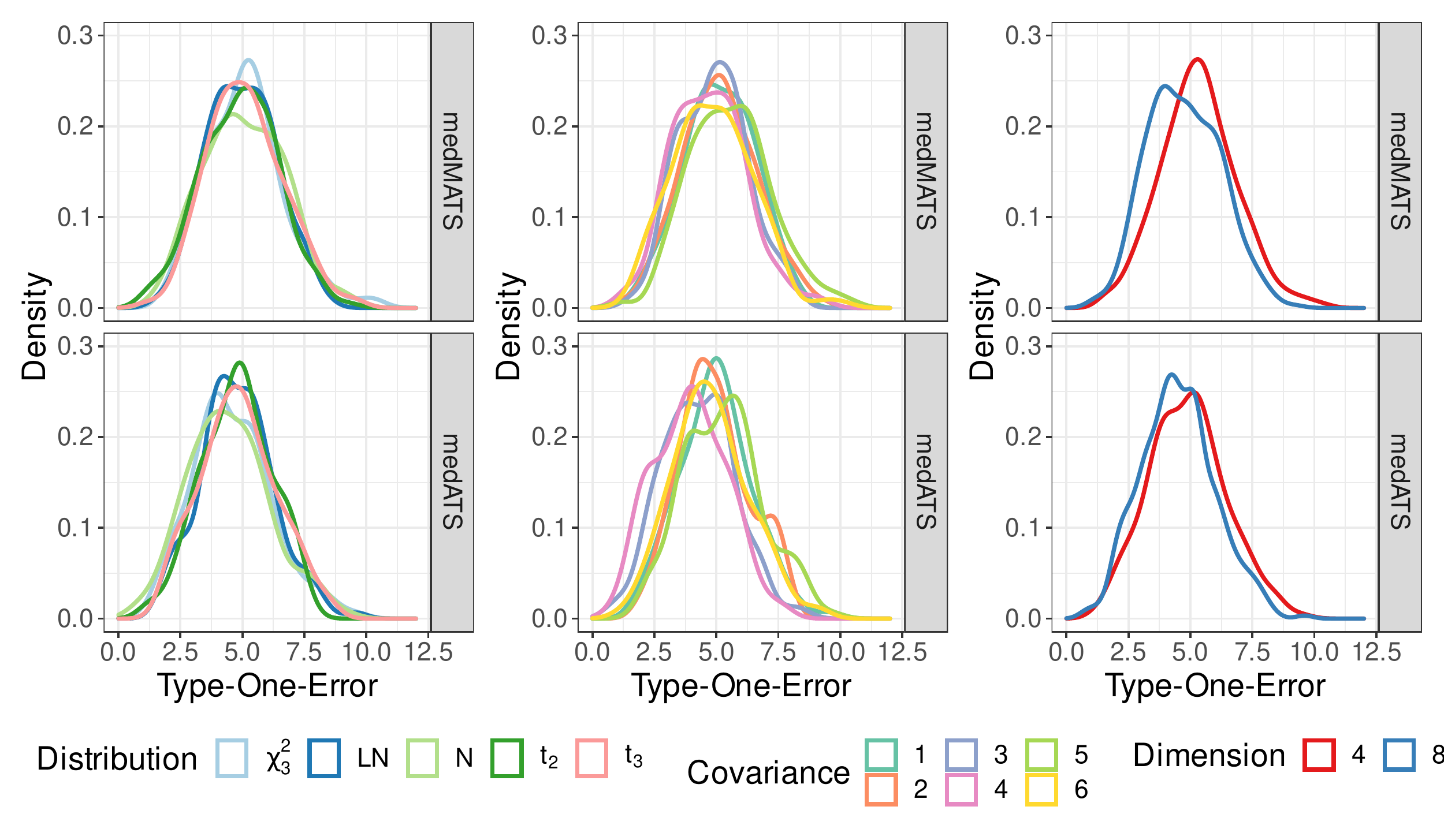}
\caption{{Kernel density estimators \citep{nadaraya_non-parametric_1965} based on Gaussian kernels} of all simulation results for the type I error of the one-way layout divided by the distribution, the covariance and the dimension of the simulation scenarios. {The bandwidth of the kernel densities is chosen by Silverman's rule-of-thumb \citep[Eq. 3.31]{silverman_density_1998}.}}
\label{graph:typeone}
\end{figure}

\begin{table}
\centering
\singlespacing
\small
\begin{tabular}{llllllllll}
  \toprule
\multirow[c]{3}{*}{$\Sigma$} & \multirow[c]{3}{*}{Distr}  & \multicolumn{4}{c}{$d=4$}  & \multicolumn{4}{c}{$d=8$}\\ \cmidrule(l{0.5em}r{0.5em}){3-6} \cmidrule(l{0.5em}r{0.5em}){7-10}
& &  \multicolumn{2}{c}{median}  & \multicolumn{2}{c}{meanMATS}& \multicolumn{2}{c}{median}  & \multicolumn{2}{c}{meanMATS}\\ \cmidrule(l{0.5em}r{0.5em}){3-4} \cmidrule(l{0.5em}r{0.5em}){5-6} \cmidrule(l{0.5em}r{0.5em}){7-8} \cmidrule(l{0.5em}r{0.5em}){9-10}
 & &  MATS & ATS & param & wild & MATS & ATS & param & wild\\
  \midrule
 \multirow[t]{6}{*}{1} & $N_{0,1}$  & 6.8 & 6.0 & $\mathbf{5.6}$ & 6.0 & 3.6 & $\mathbf{4.4}$ & $\mathbf{4.8}$ & 6.8 \\ 
 & $t_2$  & 4.0 & 3.2 & 2.8 & $\mathbf{4.8}$ & 2.8 & 3.2 & 6.0 & 6.4 \\ 
 & $t_3$  & 4.0 & 4.0 & $\mathbf{4.8}$ & 6.0 & 2.4 & 2.4 & $\mathbf{4.4}$ & $\mathbf{5.2}$ \\ 
  & $LN_{0,1}$  & $\mathbf{4.4}$ & 3.6 & - & - & $\mathbf{5.2}$ & $\mathbf{5.2}$ & - & - \\ 
 & $\chi^2_3$  & 6.8 & $\mathbf{5.2}$ & - & - & 3.2 & $\mathbf{5.2}$ & - & - \\ 
  \addlinespace
 \multirow[t]{6}{*}{2} & $N_{0,1}$  & 6.4 & $\mathbf{5.2}$ & 6.4 & 6.8 & $\mathbf{4.8}$ & $\mathbf{4.4}$ & 6.4 & 7.2 \\ 
 & $t_2$  & $\mathbf{5.2}$ & $\mathbf{5.2}$ & $\mathbf{4.4}$ & 6.8 & 4.0 & 3.6 & 4.0 & 6.8 \\ 
 & $t_3$  & 8.4 & 7.6 & $\mathbf{5.2}$ & 7.2 & 4.0 & 6.4 & 3.2 & $\mathbf{4.4}$ \\ 
 & $LN_{0,1}$  & $\mathbf{4.8}$ & $\mathbf{5.2}$ & - & - & $\mathbf{5.6}$ & 4.0 & - & - \\ 
  & $\chi^2_3$  & 4.0 & 3.2 & - & - & 2.8 & 3.2 & - & - \\ 
  \addlinespace
 \multirow[t]{6}{*}{3} & $N_{0,1}$  & 6.0 & 6.8 & $\mathbf{5.6}$ & 6.8 & 3.2 & 3.6 & $\mathbf{4.4}$ & $\mathbf{5.6}$ \\ 
  & $t_2$  & $\mathbf{4.8}$ & 4.0 & 2.8 & $\mathbf{5.2}$ & $\mathbf{4.8}$ & $\mathbf{5.2}$ & 4.0 & $\mathbf{5.6}$ \\ 
 & $t_3$  & $\mathbf{5.6}$ & $\mathbf{5.2}$ & $\mathbf{5.6}$ & 6.0 & 3.6 & $\mathbf{5.2}$ & $\mathbf{4.4}$ & 6.0 \\ 
 & $LN_{0,1}$  & 3.2 & 2.8 & - & - & $\mathbf{5.6}$ & 6.4 & - & - \\ 
 & $\chi^2_3$  & $\mathbf{4.8}$ & $\mathbf{4.4}$ & - & - & $\mathbf{4.4}$ & $\mathbf{5.2}$ & - & - \\ 
  \addlinespace
 \multirow[t]{6}{*}{4} & $N_{0,1}$  & $\mathbf{4.4}$ & 2.8 & $\mathbf{4.4}$ & $\mathbf{4.4}$ & 3.6 & 0.8 & 2.4 & 3.6 \\ 
 & $t_2$  & 2.8 & 3.6 & 2.4 & 2.4 & 3.2 & 3.2 & 1.2 & 6.0 \\ 
 & $t_3$  & 3.6 & 3.6 & 3.6 & 4.0 & $\mathbf{5.6}$ & 4.0 & 2.4 & 6.0 \\ 
 & $LN_{0,1}$  & 6.4 & $\mathbf{5.2}$ & - & - & 2.8 & 2.8 & - & - \\ 
 & $\chi^2_3$  & $\mathbf{4.8}$ & 4.0 & - & - & 3.6 & 4.0 & - & - \\ 
  \addlinespace
 \multirow[t]{6}{*}{5} & $N_{0,1}$  & $\mathbf{4.8}$ & $\mathbf{5.2}$ & $\mathbf{5.2}$ & 6.0 & 2.8 & 2.0 & 2.8 & 4.0 \\ 
 & $t_2$  & $\mathbf{4.8}$ & $\mathbf{4.4}$ & 3.6 & $\mathbf{4.4}$ & 2.8 & 2.4 & 3.2 & $\mathbf{5.6}$ \\ 
 & $t_3$  & $\mathbf{5.6}$ & 6.4 & 8.0 & 8.8 & $\mathbf{5.2}$ & 3.6 & 3.2 & $\mathbf{4.8}$ \\ 
 & $LN_{0,1}$  & 3.6 & 4.0 & - & - & 4.0 & 3.6 & - & - \\ 
 & $\chi^2_3$  & $\mathbf{5.6}$ & 6.0 & - & - & $\mathbf{5.2}$ & 4.0 & - & - \\ 
  \addlinespace
 \multirow[t]{6}{*}{6} & $N_{0,1}$  & $\mathbf{4.8}$ & 7.2 & 6.0 & 7.2 & 2.8 & 2.0 & 2.4 & 3.6 \\ 
 & $t_2$  & 4.0 & $\mathbf{4.8}$ & 2.4 & $\mathbf{4.8}$ & 2.4 & 3.2 & 0.8 & 3.6 \\ 
&  $t_3$  & 4.0 & $\mathbf{4.4}$ & 3.2 & 4.0 & 4.0 & 6.0 & 2.0 & 2.4 \\ 
 & $LN_{0,1}$  & $\mathbf{4.4}$ & 2.8 & - & - & 2.4 & $\mathbf{4.4}$ & - & - \\ 
 & $\chi^2_3$  & $\mathbf{4.8}$ & 4.0 & - & - & $\mathbf{4.4}$ & 3.2 & - & - \\ 
\bottomrule
\end{tabular}
   \caption{Type I error rate in $\%$ (nominal level $\alpha = 5\%$) for testing $\mathcal{H}_0: \mathbf{m}_1=\mathbf{m}_2$ in an one-way layout. We present all settings with the sample size $(10,\,20)$ of MATS and ATS combined with the bootstrap covariance estimator.
The two comparing meanMATS models are named with \textit{meanMATS-p} for a parametric bootstrap and with \textit{meanMATS-w} for the wild bootstrap approach.
The results inside the binomial interval for $\alpha$ $[4.4,5.6]$ are printed in bold.}
\label{table:tentwenty}
\end{table}

\begin{table}
\centering
\singlespacing
\small
\begin{tabular}{llllllllll}
  \toprule
\multirow[c]{3}{*}{$\Sigma$} & \multirow[c]{3}{*}{Distr}  & \multicolumn{4}{c}{$d=4$}  & \multicolumn{4}{c}{$d=8$}\\ \cmidrule(l{0.5em}r{0.5em}){3-6} \cmidrule(l{0.5em}r{0.5em}){7-10}
& &  \multicolumn{2}{c}{median}  & \multicolumn{2}{c}{meanMATS}& \multicolumn{2}{c}{median}  & \multicolumn{2}{c}{meanMATS}\\ \cmidrule(l{0.5em}r{0.5em}){3-4} \cmidrule(l{0.5em}r{0.5em}){5-6} \cmidrule(l{0.5em}r{0.5em}){7-8} \cmidrule(l{0.5em}r{0.5em}){9-10}
 & &  MATS & ATS & param & wild & MATS & ATS & param & wild\\
  \midrule
\multirow[t]{6}{*}{1} & $N_{0,1}$  & $\mathbf{4.8}$ & $\mathbf{4.4}$ & 6.0 & $\mathbf{5.6}$ & $\mathbf{4.8}$ & $\mathbf{4.8}$ & $\mathbf{5.2}$ & 6.0 \\ 
 & $t_2$  & 3.6 & 3.2 & $\mathbf{5.2}$ & 6.4 & 4.0 & 4.0 & 3.2 & 4.0 \\ 
 &  $t_3$  & 6.4 & $\mathbf{4.8}$ & 6.0 & 6.4 & $\mathbf{4.4}$ & $\mathbf{4.8}$ & 3.2 & 3.6 \\ 
 &  $LN_{0,1}$  & 4.0 & 3.6 & - & - & $\mathbf{4.8}$ & 6.8 & - & - \\ 
 & $\chi^2_3$  & $\mathbf{5.2}$ & 3.2 & - & - & $\mathbf{4.4}$ & $\mathbf{4.8}$ & - & - \\ 
  \addlinespace
 \multirow[t]{6}{*}{2} & $N_{0,1}$  & 3.2 & 3.6 & $\mathbf{4.4}$ & $\mathbf{4.8}$ & 3.2 & 2.8 & $\mathbf{4.4}$ & $\mathbf{4.8}$ \\ 
 & $t_2$  & 3.2 & 2.4 & 4.0 & $\mathbf{5.6}$ & 2.0 & 2.4 & 1.6 & $\mathbf{4.8}$ \\ 
 & $t_3$  & 4.0 & $\mathbf{5.2}$ & 2.4 & 2.8 & 2.4 & 2.0 & 3.2 & 4.0 \\ 
 & $LN_{0,1}$  & $\mathbf{5.2}$ & $\mathbf{5.6}$ & - & - & 4.0 & $\mathbf{4.4}$ & - & - \\ 
 & $\chi^2_3$  & 7.6 & 7.2 & - & - & 4.0 & 3.2 & - & - \\ 
  \addlinespace
 \multirow[t]{6}{*}{3} & $N_{0,1}$  & 2.8 & 2.0 & 2.8 & 2.8 & $\mathbf{4.8}$ & 3.2 & 3.6 & 6.4 \\ 
 & $t_2$  & $\mathbf{5.6}$ & 6.4 & 3.6 & 6.8 & 3.2 & 2.0 & 1.6 & 2.4 \\ 
 & $t_3$  & $\mathbf{5.2}$ & $\mathbf{4.4}$ & $\mathbf{4.4}$ & $\mathbf{5.6}$ & 3.2 & 3.2 & 1.6 & 4.0 \\ 
 & $LN_{0,1}$  & 4.0 & 2.4 & - & - & 3.6 & 3.2 & - & - \\ 
 & $\chi^2_3$  & 6.0 & 6.8 & - & - & 3.2 & 3.6 & - & - \\ 
  \addlinespace
 \multirow[t]{6}{*}{4} & $N_{0,1}$  & 3.6 & 2.4 & $\mathbf{4.8}$ & $\mathbf{5.6}$ & 3.6 & 1.6 & 4.0 & $\mathbf{5.6}$ \\ 
 & $t_2$  & 1.6 & 1.6 & 2.0 & 2.8 & 1.2 & 1.2 & 1.6 & 4.0 \\ 
 & $t_3$  & $\mathbf{5.6}$ & 2.4 & 1.6 & 1.6 & $\mathbf{4.4}$ & 2.0 & 2.0 & $\mathbf{5.6}$ \\ 
 & $LN_{0,1}$  & $\mathbf{4.4}$ & 2.0 & - & - & 2.8 & 2.0 & - & - \\ 
 & $\chi^2_3$  & $\mathbf{4.8}$ & 4.0 & - & - & $\mathbf{5.2}$ & 3.2 & - & - \\ 
  \addlinespace
 \multirow[t]{6}{*}{5} & $N_{0,1}$  & 8.8 & 8.4 & 7.6 & 8.4 & $\mathbf{4.8}$ & 3.6 & 4.0 & $\mathbf{4.8}$ \\ 
 & $t_2$  & 4.0 & 2.8 & 3.2 & 4.0 & 6.8 & $\mathbf{5.2}$ & 4.0 & $\mathbf{4.4}$ \\ 
 & $t_3$  & 9.6 & 8.4 & 6.8 & 7.6 & $\mathbf{4.8}$ & $\mathbf{4.4}$ & 6.4 & 7.2 \\ 
 & $LN_{0,1}$  & 6.0 & $\mathbf{5.2}$ & - & - & $\mathbf{5.2}$ & 6.0 & - & - \\ 
 & $\chi^2_3$  & 6.8 & 6.8 & - & - & $\mathbf{5.6}$ & 6.4 & - & - \\ 
  \addlinespace
 \multirow[t]{6}{*}{6} & $N_{0,1}$  & $\mathbf{5.6}$ & 6.0 & $\mathbf{5.6}$ & 6.4 & 2.4 & 2.8 & 2.0 & 3.2 \\ 
 & $t_2$  & 6.4 & 4.0 & 3.2 & 6.8 & 3.2 & 4.0 & 1.6 & 3.2 \\ 
 & $t_3$  & 6.4 & 7.2 & 4.0 & 6.8 & 4.0 & $\mathbf{4.4}$ & 2.4 & 6.8 \\ 
 & $LN_{0,1}$  & 4.0 & $\mathbf{4.8}$ & - & - & 2.0 & 3.6 & - & - \\ 
 & $\chi^2_3$  & 2.4 & 2.0 & - & - & 3.6 & 2.4 & - & - \\ 
  \bottomrule
  \end{tabular}
 \caption{Type I error rate in $\%$ (nominal level $\alpha = 5\%$) for testing $\mathcal{H}_0: \mathbf{m}_1=\mathbf{m}_2$ in an one-way layout. We present all settings with the sample size $(20,\,10)$ of MATS and ATS combined with the bootstrap covariance estimator.
The two comparing meanMATS models are named with \textit{meanMATS-p} for a parametric bootstrap and with \textit{meanMATS-w} for the wild bootstrap approach.
The results inside the binomial interval for $\alpha$ $[4.4,5.6]$ are printed in bold.}
\label{table:twentyten}
\end{table}

To get a deeper insight into the tests performance, we present the result of the type I error simulation study as density plots divided by certain aspects.
The left column in Figure \ref{graph:typeone} presents the density divided by the different distributions, the middle column by the covariance settings and the right column by the dimension.
Thus, it becomes clear which aspects influence the tests performance.
The right column shows that the different dimensions have almost none influence on the tests' performance in both methods.
On the contrary, there are clear differences between the covariance settings. 
They seem to have the most influence on the type I errors in this comparison.
Here, the type I errors for the ATS statistic with bootstrap covariance estimator differs more than the errors for the MATS statistic.
This phenomenon can also be observed in the first column of Figure \ref{graph:typeone}; and even for the standard normal distribution, the ATS statistic produces in average too small type  I errors.
This can also be observed in Figure \ref{graph:boxplots} and \ref{graph:comp} of the paper.
All simulation results of the type  I error for the sample sizes $(10,\;20)$ and $(20,\;10)$ are presented in in Tables \ref{table:tentwenty} and \ref{table:twentyten}.
{There, it can be seen, that the unbalanced sample size seems to have no influence of the tests performance.}
Tendentiously, the behaviour gets better as the sample size gets larger.
All in all, Figure \ref{graph:typeone} shows that the QMANOVA method is still robust against various distributional scenarios. 
This is consistent with our theoretical results.

\subsection{Two-way-Layout}
\begin{figure}
\includegraphics[scale=0.56]{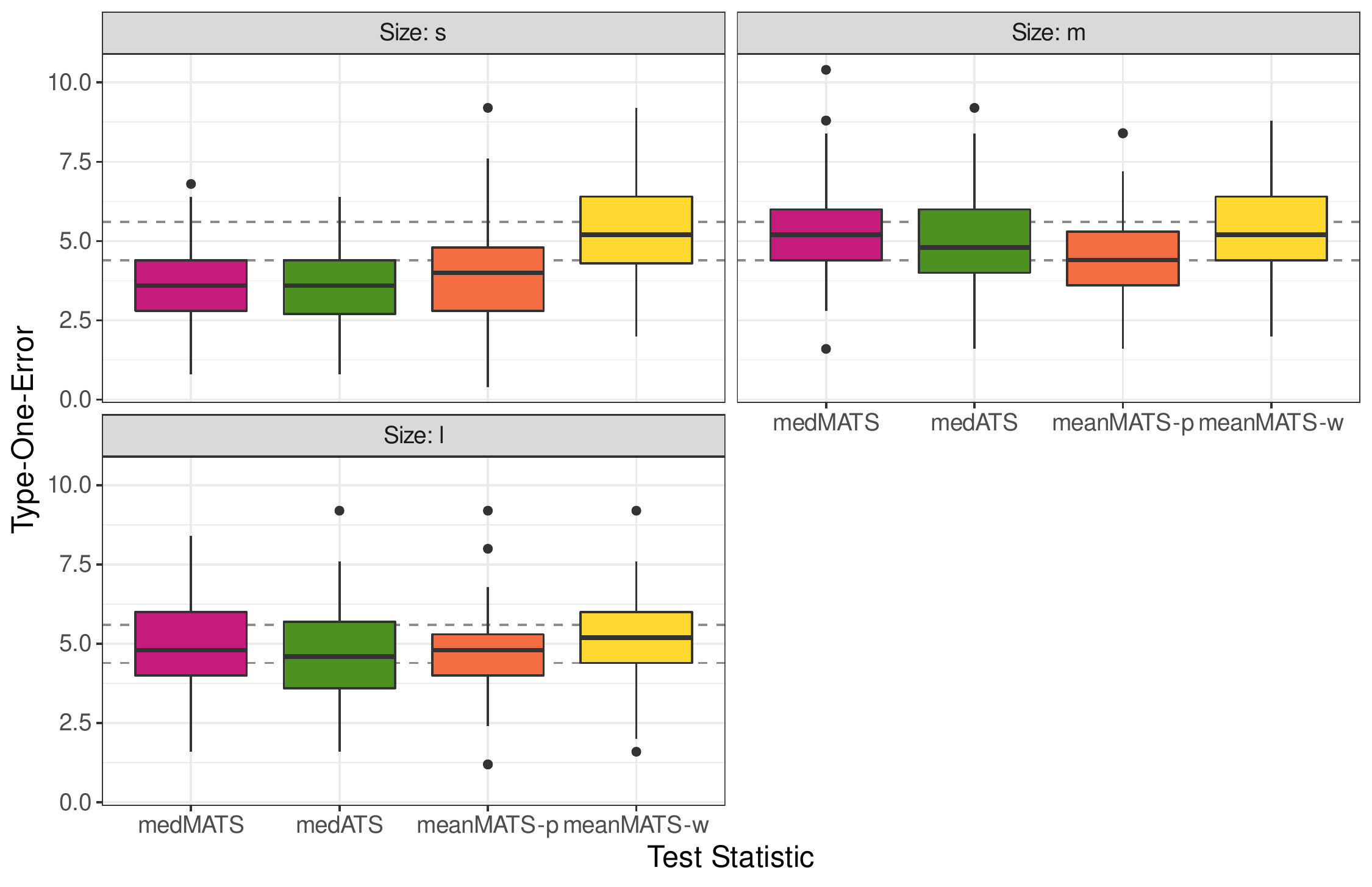}
\centering
\caption{Type I error (in percent) for the symmetric distributions ($N_{0,1}$, $t_2$, $t_3$)  in the two-way-layout simulations study for the QMANOVA with bootstrap covariance estimator and the mean-based MATS by \citet{friedrich_mats_2018} with parametric (p) and wild (w) bootstrap as resampling technique. {The results are divided by sample size. Here, the samples $n^{(4)},\,n^{(5)}$ and $n^{(6)}$ are categorized as small (\textit{s}), the samples $n^{(7)}$ and $n^{(10)}-n^{(12)}$ as middle (\textit{m}) and the samples $n^{(8)}$ and $n^{(9)}$ as large (\textit{l}).}}
\label{graph:twowaycomp}
\end{figure}

\begin{figure}
\includegraphics[scale=0.7]{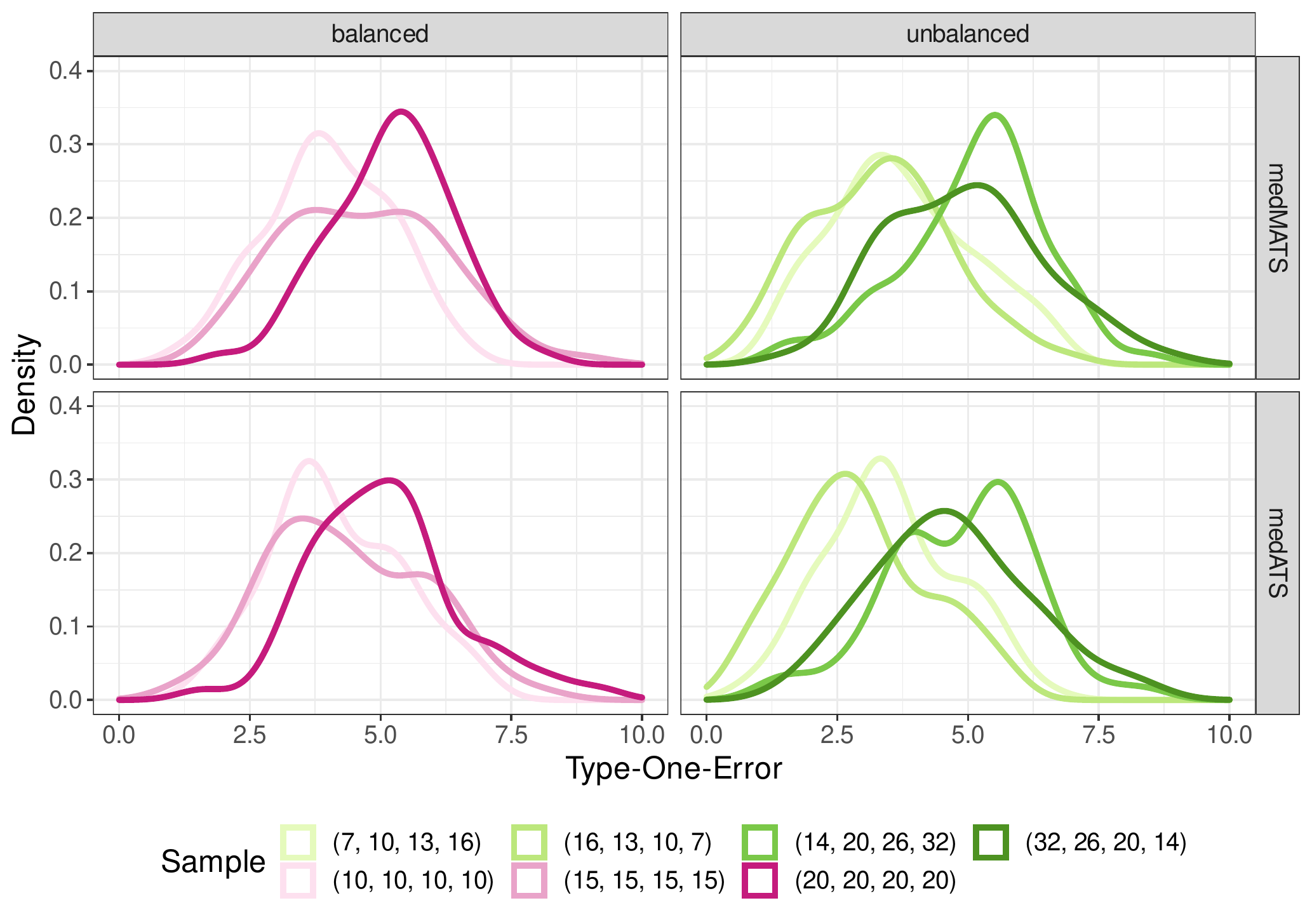}
\caption{{Kernel density estimators \citep{nadaraya_non-parametric_1965} based on Gaussian kernels} of all simulation results for the type I error of the two-way layout divided by the sample sizes. {The bandwidth of the kernel densities is chosen by Silverman's rule-of-thumb \citep[Eq. 3.31]{silverman_density_1998}.}}
\label{graph:sample}
\end{figure}

To get more insights into our method, we analyse the behaviour in the case of a two-way-layout in the following set-up.
We consider two crossed factors $A$ and $B$ having two levels each with the sample sizes {$n^{(4)}=(n_{ab})_{a,b=1,2}=(n_{11},n_{12}, n_{21},n_{22})=(7, 10, 13, 16),$} $\;n^{(5)}=(10, 10, 10, 10)$, $n^{(6)}=(16, 13, 10, 7),$ $n^{(7)}=(20, 20, 20, 20),$ $n^{(8)}=(100, 100, 100, 100)$, $n^{(9)}=(100, 50, 100, 50)$, $n^{(10)}=(15, 15, $ $15, 15)$, $n^{(11)}=(14, 20, $ $16, 32)$ and $n^{(12)}=(32, 26, 20, 14)$.
The data generating process and the distributions are identical to the simulation study for a one-way layout and use the following covariance settings with $d=4$ dimensions:
\begin{enumerate}
\setcounter{enumi}{6}
	\item $\mathbf{\Sigma}^{(i_1i_2)}=\mathbf{I}_d-\frac{1}{2}\left(\mathbf{J}_d-\mathbf{I}_d\right),\;i_1,i_2\in\{1,2\}$,
	\item $\mathbf{\Sigma}^{(i_1i_2)}=\left(0.6^{|a-b|}\right)_{a,b=1}^d,\;i_1,i_2\in\{1,2\}$,
	\item $\mathbf{\Sigma}^{(i_1i_2)}=\xi\,\mathbf{I}_d+\frac{1}{2}\left(\mathbf{J}_d-\mathbf{I}_d\right),\;\xi\in\{1,\dots,4\}$,
	\item $\mathbf{\Sigma}^{(i_1i_2)}=\left(0.6^{|a-b|}\right)_{a,b=1}^d+\xi\,\mathbf{I}_d,\;\xi\in\{1,\dots,4\}.$
	\end{enumerate}
We consider three different hypotheses to the significance level $\alpha=0.05$:
\begin{itemize}
\item The hypothesis of no effect of factor $A$ 
\[\mathcal{H}_0(A):\left(\mathbf{P}_2\otimes\frac{1}{2}\mathbf{J}_2\otimes\mathbf{I}_4\right)\mathbf{q}=\mathbf{0}_{8}\Leftrightarrow\mathcal{H}_0(A):\mathbf{\bar q}_{1\cdot}=\mathbf{\bar q}_{2\cdot},\]
\item the hypothesis of no effect of factor $B$
\[\mathcal{H}_0(B):\left(\frac{1}{2}\mathbf{J}_2\otimes\mathbf{P}_b\otimes\mathbf{I}_4\right)\mathbf{q}=\mathbf{0}_{8}\Leftrightarrow\mathcal{H}_0(B):\mathbf{\bar q}_{\cdot 1}=\mathbf{\bar q}_{\cdot 2},\]
\item the hypothesis of no $A\times B$ interaction effect
\[\mathcal{H}_0(AB):\left(\mathbf{P}_2\otimes\mathbf{P}_2\otimes\mathbf{I}_4\right)\mathbf{q}=\mathbf{0}_{8}.\]
\end{itemize}
This leads to $540$ different scenarios and the whole set-up is similar to the simulations in \citet{friedrich_mats_2018}; the mean-based MATS is our comparing method again.
In general, QMANOVA's performance for the two-way layout is similar to the performance in the one-way-layout.
Again, the bootstrap covariance estimator performs best and the performance of the two test statistics is similar.
In the upper graph of Figure \ref{graph:twowaycomp} this result is shown for the symmetric distributions.
The performance of medMATS is similar to the meanMATS-w and analogously to the simulations regarding the one-way layout, the medMATS tends to be conservative and the meanMATS-w tends to be more liberal.
With the density plots in the lower graph of Figure \ref{graph:twowaycomp}, we analyse again the influence of the distribution, the covariance and the test layout on the test's performance.
In the distribution and the covariance column, one can see more variation than in the equivalent figures of the one-way layout.
{Since the layout is more complex, this observation is not surprising and leads us to the question if the two-way layout works with the very small sample size of $(10,\,10)$ as the one-way layout does.}
{In the situation of such small samples}, the QMANOVA methods tend to be conservative.
Figure \ref{graph:sample} compares the tests performance regarding the sample size.
Here, the largest sample sizes $n^{(8)}=(100, 100, 100, 100)$ and $n^{(9)}=(100, 50, 100, 50)$ are excluded from the analyses.
It is apparent from this figure, that the tests performance is very well for the sample sizes $n^{(7)}=(20, 20, 20, 20),\;n^{(11)}=(14, 20, 16, 32)$ and $n^{(12)}=(32, 26, 20, 14)$.
The smaller samples seems to be too small for the QMANOVA, because the simulations study produces there fairly too conservative type I errors.
For the two-way layout, most empirical type I errors under $4.4$, the lower bound of the $95\%$ binomial interval, are from small sample sizes.
{Furthermore, it can be observed, that the QMANOVA can handle unbalanced scenarios equally well as the balanced ones.}
As in the settings regarding a one-way layout, the medMATS with the bootstrap covariance estimator is the best QMANOVA method.
Here, the settings with conservative behaviour (empirical type I error smaller than $4.4$) have often covariance settings $9$ or $10$, if they have not small sample sizes.
Again, conservative values appear more frequently than liberal one (empirical type I error bigger than $5.6$).
Looking at all results with liberal behaviour, there is the reversed situation as for the conservative results: the settings with liberal behaviour have often covariance settings $7$ and $8$.
From this it follows that unequal covariances seem to be a problem in a two-way layout, this was not observable in the one-way layout.
{In summary, the unequal covariances have a clear effect on the simulation results in the two-way layout, which was not the case in the one-way layout.}

\begin{figure}
\includegraphics[scale=0.57]{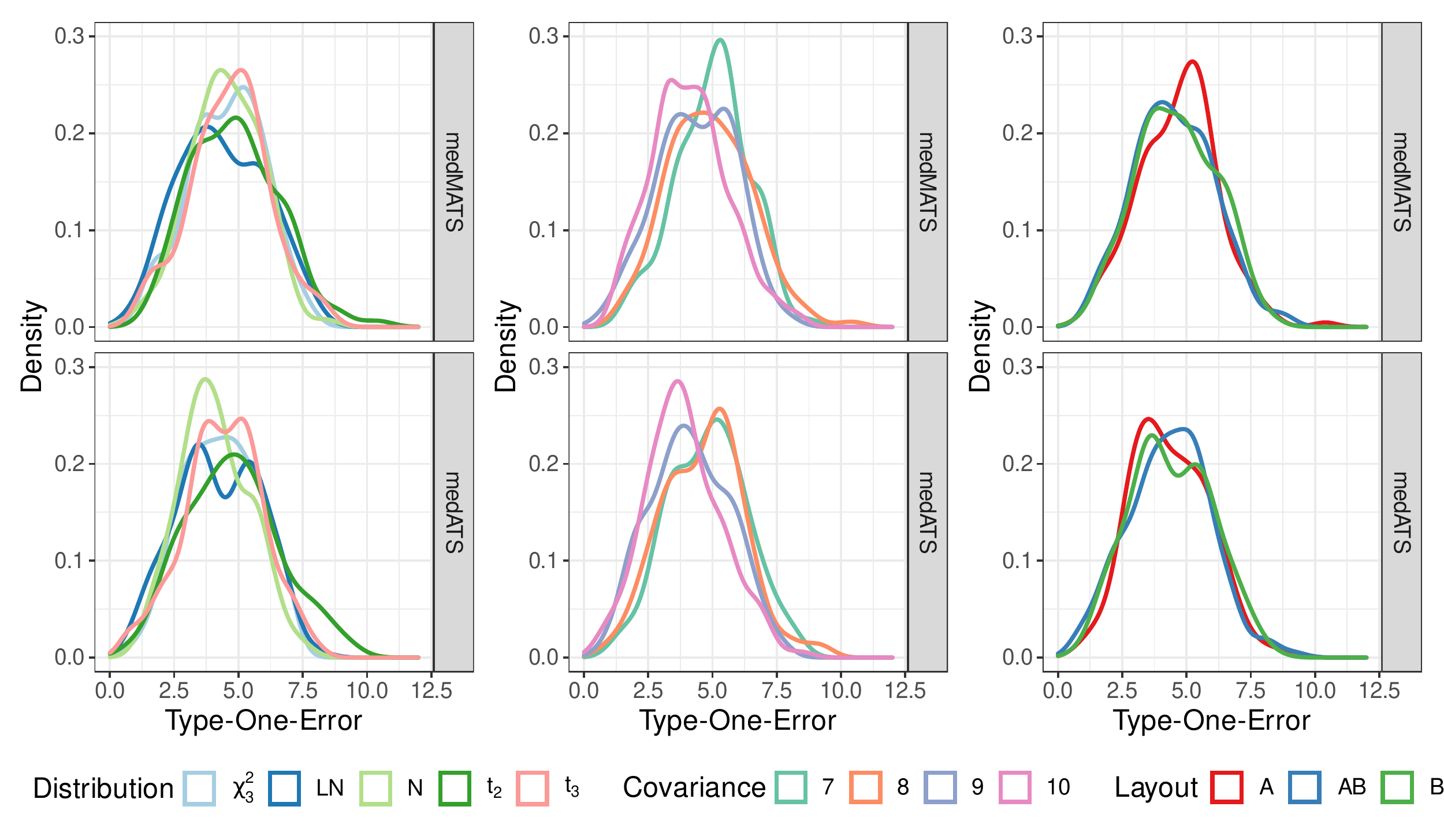}
\centering
\caption{{Kernel density estimators \citep{nadaraya_non-parametric_1965} based on Gaussian kernels} of all simulation results for the type I error of the two-way layout divided by the distribution, the covariance and the layouts of the simulation scenarios. {The bandwidth of the kernel densities is chosen by Silverman's rule-of-thumb \citep[Eq. 3.31]{silverman_density_1998}.}}
\label{graph:twowaydiff}
\end{figure}

\section{Mathematical Foundation of the Proposition \ref{jointconsistent}}\label{mathfound}
The vector of group-specific marginal empirical distribution functions is equal to
	\begin{align}\label{empmeasure}
	\begin{pmatrix}
	\hat{F}_{i1}(t_1)\\ \vdots\\ \hat{F}_{id}(t_d)
	\end{pmatrix}
	=\frac{1}{n_i}\sum_{j=1}^{n_i}\begin{pmatrix}
	\mathbf{1}_{(-\infty,t_1]}(X_{ij1})\\ \vdots\\ \mathbf{1}_{(-\infty,t_d]}(X_{ijd})
	\end{pmatrix}
	=\frac{1}{n_i}\sum_{j=1}^{n_i}\begin{pmatrix}
	\mathbf{1}_{(-\infty,t_1]\times\br^{d-1}}(X_{ij1},\dots,X_{ijd})\\ \vdots\\ \mathbf{1}_{\br^{d-1}\times(-\infty,t_d]}(X_{ij1},\dots,X_{ijd})
	\end{pmatrix}.
	\end{align}
 For $t\in\br$ and $\ell,r\in\{1,\dots,d\}$, set $A_{t\ell r}=(-\infty,t]$ if $\ell=r$ and $A_{t\ell r}=\br$ else.
Then, the functions
	\begin{align*}
	f_{t\ell}:=\mathbf{1}_{\bigotimes_{r=1}^d A_{t\ell r}}:\br^d\to\br,\,
	\begin{pmatrix}x_1\\ \vdots\\ x_d\end{pmatrix}\mapsto f_{t\ell}(x_1,\dots,x_d)&=\mathbf{1}_{\bigotimes_{r=1}^d A_{t\ell r}}(x_1,\dots,x_d)\\
	&=\begin{cases}
	1,\;x_\ell\in A_{t\ell\ell}=(-\infty,t]\\
	0,\;x_\ell\not\in A_{t\ell\ell}=(-\infty,t]
	\end{cases}
	\end{align*}
characterize the marginal empirical distribution functions and they form a set:
	\begin{align*}
	\cale:=\left\{f_{t\ell}=\mathbf{1}_{\bigotimes_{r=1}^d A_{t\ell r}}| A_{t\ell\ell}=(-\infty,t]\land A_{t\ell r}=\br:r\not=\ell;\,r,\ell\in\{1,\dots,d\}\land t\in\br\right\}.
	\end{align*}

\begin{lemma}\label{VC}
The set of measurable functions $\mathcal{E}$ is a VC-class.
\end{lemma}

\begin{proof}
	Due to Problem 9 in Section 2.6 of \citet[p.~151]{vaart_weak_2000} it remains to show that 
		\[\mathcal{C} = \left\{ \bigotimes_{r=1}^d A_{t\ell r}| A_{t\ell\ell}=(-\infty,t] \land A_{t\ell r}=\br:\,r\neq \ell\in\{1,\ldots,d\}\land t\in\br\right\}\]
		is a VC-class.
	We use the fact that the $d$-dimensional cells 
	\[\mathcal D = \left\{\bigotimes_{r=1}^d A_r| A_r=(-\infty,t_r]:r\in\{1,\ldots,d\}\land t_r\in\br\right\}\]
	 form a VC-class with VC-index $d+1$ \citep[cf.][Ex. 2.6.1]{vaart_weak_2000}. 
	To prove that $\mathcal{C}$ is a VC-class, let us suppose for a moment that $\mathcal{C}$ shatters the subset $\{\mathbf{x}_1,\ldots,\mathbf{x}_m\}\subset\br^d$ for some $m\in\bn$, e.g. every subset $\mathbf{X}$ of $\{\mathbf{x}_1,\ldots,\mathbf{x}_m\}$ can be written as an intersection between $\{\mathbf{x}_1,\ldots,\mathbf{x}_m\}$ and an element $C\in\mathcal{C}$: $\mathbf{X}=\{\mathbf{x}_1,\ldots,\mathbf{x}_m\}\cap C$.
	Define $M= \max\left\{x_{sr}| s\in\{1,\ldots,m\};r\in\{1,\dots,d\}\right\} + 1$.
	Thus, $M$ is larger than any component of $\mathbf{x}_1,\ldots,\mathbf{x}_m$.
	Then it is clear that 
		\[\mathcal{C}'=\left\{\bigotimes_{r=1}^d A_{t\ell r}| A_{t\ell\ell}=(-\infty,t]\land A_{tr\ell}=(-\infty,M]:r\neq\ell\in\{1,\ldots,d\}\land t\in\br\right\}\]
	shatters $\{\mathbf{x}_1,\ldots,\mathbf{x}_m\}\subset \br^d$ as well. 
	Since $\mathcal C'$ is clearly a subset of $\mathcal D$, $\mathcal{C}'$ is for $m\le d+1$ also a VC-class with VC-index $d+1$ or smaller.
	This is a contradiction.
	Thus, $\mathcal{C}$ does not shatter the subset $\{\mathbf{x}_1,\ldots,\mathbf{x}_m\}\subset\br^d$ and it can finally be deduced that $\mathcal C$ is a VC-class with a VC-index smaller or equal to $d+1$.
	\end{proof}
The set $\cale$  corresponds to the marginal empirical distribution functions because the empirical measure of the functions in $\cale$ create them, compare Formula (\ref{empmeasure}). 
Categorizes the set $\cale$ as a VC-Class amounts to apply the theory of empirical processes on the marginal empirical distribution functions.
Applying Theorem 2.6.7 in \citet{vaart_weak_2000}, Lemma \ref{VC} yields that $\mathcal{E}$ satisfies the uniform entropy condition (2.5.1) in \citet[p.~127]{vaart_weak_2000}.
The function
\[\mathbf{1}_{\br^ d}:\br^d\to\br,\;(x_1,\dots,x_d)\mapsto 1\]
is an envelope function for $\mathcal{E}$ and it holds for every distribution $P$ on $\left(\br^d,\mathcal{B}\left(\br^d\right)\right)$
	\begin{align*}
	\left(\int|\mathbf{1}_{\br^d}|^2dP\right)^{\frac{1}{2}}=\left(\int_{\br^d}\mathbf{1}dP\right)^{\frac{1}{2}}=(P(\br^d))^{\frac{1}{2}}=1<\infty.
	\end{align*}
From these three conditions it can be concluded that $\mathcal{E}$ is $P$-Donsker \citep[cf.][p.~141]{vaart_weak_2000}.
Consequently:
	\begin{align*}
	\sqrt{n_i}\left(\hat{F}_{i\ell}-F_{i\ell}\right)_{\ell\in\{1,\dots,d\}}\stackrel{d}{\longrightarrow}\mathbb{G}_i\;\text{in}\;\ell^\infty(\cale),\;i\in\{1,\dots,k\}.
	\end{align*}
By construction the limit process $\mathbb{G}_i:\mathcal{E}\to\br,\,f\mapsto \mathbb{G}_if$ is a empirical process in $\ell^\infty(\mathcal{E})$ for every $i\in\{1,\dots,k\}$.
It can also be written as $\{\mathbb{G}_if_{t\ell}|f_{t\ell}\in\cale\}$  and is furthermore a zero-mean Gaussian process \citep[cf.][p.~81 f.]{vaart_weak_2000}.
Its covariance is constructed as follows:
Let $f_{t\ell}=\mathbf{1}_{\bigotimes_{r=1}^dA_{t\ell r}},\,f_{sm}=\mathbf{1}_{\bigotimes_{r=1}^dB_{smr}}\in\cale$. 
Thus, it holds $A_{t\ell\ell}=(-\infty,t],\,B_{smm}=(-\infty,s]$ and $A_{t\ell r}=B_{sm r}=\br$ for the other $r\in\{1,\dots,d\}$ and \citep[cf.][Equation (2.1.2)]{vaart_weak_2000}
	\begin{align*}
	\cov(\mathbb{G}_if_{t\ell},\mathbb{G}_if_{sm})&=E\mathbb{G}_if_{t\ell}\mathbb{G}_if_{sm}=P_if_{t\ell}f_{sm}-P_if_{t\ell}P_if_{sm}\\
	&=P_i\left(\bigotimes_{r=1}^dA_{t\ell r}\cap\bigotimes_{r=1}^dB_{smr}\right)-P_i\left(\bigotimes_{r=1}^dA_{t\ell r}\right)P_i\left(\bigotimes_{r=1}^dB_{smr}\right)\\
	&=P_i\left(\bigotimes_{r=1}^d\left(A_{t\ell r}\cap B_{smr}\right)\right)-P_i\left(\bigotimes_{r=1}^dA_{t\ell r}\right)P_i\left(\bigotimes_{r=1}^dB_{smr}\right).
	\end{align*}
For $\ell=m$ it follows that $A_{t\ell\ell}=(-\infty,t]$ and $B_{smm}=(-\infty,s]$ and thus $A_{t\ell\ell}\cap B_{smm}=(-\infty,\min(t,s)]$.
And for $\ell\not=m$ follows that $A_{t\ell\ell}=(-\infty,t],\,B_{sm\ell}=\br,\,A_{t\ell m}=\br$ and $B_{smm}=(-\infty,s]$ and thus $A_{t\ell\ell}\cap B_{sm\ell}=(-\infty,t]$ and $A_{t\ell m}\cap B_{smm}=(-\infty,s]$.
For the covariance this means:
	\begin{align*}
	\cov(\mathbb{G}_if_{t\ell},\mathbb{G}_if_{sm})=
	\begin{cases}
	F_{i\ell}\left(\min(t,s)\right)-F_{i\ell}(s)F_{i\ell}(t),\;\ell=m\\
	F_{im\ell}(s,t)-F_{i\ell}(t)F_{im}(s),\;\ell\not=m
	\end{cases}.
	\end{align*}
By Lemma 1.5.3 in \citet{vaart_weak_2000} this characterizes $\mathbb{G}_i$ in $\ell^\infty(\cale)$ completely.
For a fixed $\ell\in\{1,\dots,d\}$ any function $f_{t\ell}\in\cale$ can be identified with $t\in\br$, and thus, $\ell^\infty(\cale)$ with $\left[\ell^\infty\left(\br\right)\right]^d$ \citep[cf.][Example 2.1.3]{vaart_weak_2000}.
The Skorokhod Space $D(\br)$ contains all right-continuous functions $G:\br\to\br$ with left limits \citep[cf.][p.~3]{vaart_weak_2000}. 
When the space is equipped with the supremum norm $\| G\|=\sup_{t\in\br}|G(t)|$, the weak convergence in $D(\br)^d$ follows from the weak convergence in $\left[\ell^\infty\left(\br\right)\right]^d$ \citep[cf.][Example 2.1.3]{vaart_weak_2000}.
Thus, the convergence in distribution holds also in $D(\br)^d$:
	\begin{align}\label{multidonsker}
	\sqrt{n_i}\left(\hat{F}_{i\ell}-F_{i\ell}\right)_{\ell\in\{1,\dots,d\}}\stackrel{d}{\longrightarrow}\mathbb{G}_i\;\text{in}\;D(\br)^d,\;i\in\{1,\dots,k\}.
	\end{align}
\section{Consistency of the Versions of $E(\mathbf{T},\mathbf{\hat\Sigma})$}
\begin{prop}\label{consistent}
	For a consistent estimator $\mathbf{\hat \Sigma}$ of $\mathbf{\Sigma}$ the versions of $E(\mathbf{T},\mathbf{\hat\Sigma})$ in the ATS and the MATS are consistent for $E\left(\mathbf{T},\mathbf{\Sigma}\right)$, thus,
	\begin{enumerate}
		\item $\tr\left(\mathbf{T\hat\Sigma T}\right)^{-1}\mathbf{I}_{dk}\stackrel{P}{\longrightarrow}\tr\left(\mathbf{T\Sigma T}\right)^{-1}\mathbf{I}_{dk}$;
		\item $\left(\mathbf{T}\mathbf{\hat\Sigma}_0\mathbf{T}\right)^+\stackrel{P}{\longrightarrow}\left(\mathbf{T}\mathbf{\hat\Sigma}_0\mathbf{T}\right)^+$.
	\end{enumerate}			 
\end{prop}
\begin{proof}
	The estimator $E(\mathbf{T},\mathbf{\hat\Sigma})=\tr(\mathbf{T\hat\Sigma T})^{-1}\mathbf{I}_{dk}$ is consistent for $\tr\left(\mathbf{T\Sigma T}\right)^{-1}\mathbf{I}_{dk}$ as a continuous function of the consistent estimator $\mathbf{\hat\Sigma}$ for $\mathbf{\Sigma}$. 
	Instead of the classical inverse, the Moore-Penrose inverse is not a continuous function.
	That is why there is more to do to prove the consistency of $E(\mathbf{T},\mathbf{\hat\Sigma})=(\mathbf{T}\mathbf{\hat\Sigma}_0\mathbf{T})^+$.
	The consistency of $\mathbf{\hat\Sigma_0}$ follows from the consistency of $\mathbf{\hat\Sigma}$.
	And by the Continuous Mapping Theorem \citep[Thm. 1.11.1]{vaart_weak_2000} it holds $\mathbf{T}\mathbf{\hat\Sigma}_0\mathbf{T}\stackrel{P}{\to}\mathbf{T}\mathbf{\Sigma}_0\mathbf{T}$ and $\mathbf{\hat\Sigma_0}$ has full rank.
	Consequently, there is no rank jump in $\mathbf{T}\mathbf{\hat\Sigma}_0\mathbf{T}$ and from Theorem 4.2 in \citet{rakocevic_continuity_1997} it follows the assertion.
\end{proof}
\section[Bootstrap Versions of the Covariance Estimators]{Bootstrap Versions of the Covariance Estimators}\label{bootstrapversions}
	To verify the consistency of the bootstrap versions of the three covariance estimators, one has to verify the consistency of the different plug-ins.
		\begin{lemma}\label{bootjointconsistent}
		Under Assumption \ref{jointcontinuous} the estimator $\hat{F}^*_{i\ell m}\left(\hat{q}^*_{i\ell},\hat{q}^*_{im}\right)$ is given the data strong 		consistent for $F_{i\ell m}\left(q_{i\ell},q_{im}\right)$ for every $i\in\{1,\dots,k\}$, $\ell,m\in\{1,\dots,d\}$.
		\end{lemma}
		\begin{proof}
	The statement is proved analogously to Proposition \ref{jointconsistent}.
	By Example 2.1.3 of \citet{vaart_weak_2000} the set $\{\mathbf{1}_{(-\infty,t_\ell]\times(-\infty,t_m]}|(t_\ell,t_m)\in\br^2\}$ forms a Donsker-Class.
	As a consequence, with Theorem 3.6.2 of \citet{vaart_weak_2000} given the data in probability applies:
		\begin{align*}
		\sqrt{n_i}\left(\hat{F}_{i\ell m}^*-\hat{F}_{i\ell m}\right)_{\ell,m\in\{1,\dots,d\}}\stackrel{d}{\longrightarrow}\mathbb{G}_i\;\text{in}\;D(\br)^d,\;i\in\{1,\dots,k\}.
		\end{align*}
	Analogous to the proof of Proposition \ref{jointconsistent} it follows
		\begin{align*}
		\sup_{(t_1,t_2)\in\br^2}\left|\hat F^*_{i\ell m}(t_1,t_2)- \hat F_{i\ell m}(t_1,t_2)\right|\stackrel{a.e.}{\longrightarrow}0
		\end{align*}
	and from this
		\[\sup_{(t_1,t_2)\in\br^2}\left|\hat F^*_{i\ell m}(t_1,t_2)- F_{i\ell m}(t_1,t_2)\right|\stackrel{a.e.}{\longrightarrow}0.\]
	The assertion ensues with the same arguments as in the rest of the Proposition \ref{jointconsistent} proof.
	\end{proof}
	Analogous to Lemma S.2 in the supplement of \citet[p.~12]{ditzhaus_qanova_2021} one can prove the consistency of the bootstrap counterpart of the interval-based estimator:
	\begin{lemma}\label{bootinterval-based}
	Let $\hat{\sigma}_{i\ell}^{PB*}(p)$ be the bootstrap counterpart of the interval-based estimator $\hat\sigma_{i\ell}^{PB}(p)$ defined in (\ref{interval-based}).
	Then, it holds the following conditional convergence given the data in probability:
		\begin{align*}
		\hat{f}_{i\ell}^*=\frac{\sqrt{p(1-p)}}{\hat{\sigma}_{i\ell}^{PB*}(p)}\stackrel{P}{\longrightarrow}f_{i\ell}\left(F_{i\ell}^{-1}(p)\right),\, i\in\{1,\dots,k\},\ell\in\{1,\dots,d\}.
		\end{align*}
	\end{lemma}
	\begin{proof}
	Let the observations be fixed. 
	Having subsequence arguments in mind, we can assume without loss of generality that \eqref{bootmultidonsker} holds.
	For a fixed $i\in\{1,\dots,k\}$ and $\ell\in\{1,\dots,d\}$ recognize that
		\begin{align}
		\hat\sigma_{i\ell}^{PB*}(p)&=\sqrt{n_i}\;\cfrac{X_{u_i(p):n_i}^{(i\ell)*}-X_{l_i(p):n_i}^{(i\ell)*}}{2z_{\alpha_{n_i\ell}^*(p)/2}+2n_i^{-1/2}}\nonumber\\
		\Leftrightarrow\left(2z_{\alpha_{n_i\ell}^*(p)/2}+2n_i^{-1/2}\right)\hat\sigma_{i\ell}^{PB*}(p)&=\sqrt{n_i}\left(X_{u_i(p):n_i}^{(i\ell)*}-X_{l_i(p):n_i}^{(i\ell)*}\right)\nonumber\\
		&=\sqrt{n_i}\left(\hat{F}_{i\ell}^{*-1}\left(\frac{u_i(p)}{n_i}\right)-\hat{F}_{i\ell}^{*-1}\left(\frac{l_i(p)}{n_i}\right)\right)\nonumber\\
		&=\sqrt{n_i}\left(\hat{F}_{i\ell}^{*-1}\left(p+k_{n_i,1}\right)-\hat{F}_{i\ell}^{*-1}\left(p-k_{n_i,2}\right)\right),\label{formularbootint}
		\end{align}
	where $\sqrt{n_i}k_{n_i,j}\to\xi:=z_{\frac{\alpha}{2}}\sqrt{p(1-p)},\,j\in\{1,2\}$.
	This approach is taken from the proof of Lemma S.2 in the supplement of \citet{ditzhaus_qanova_2021}.
	By the definition of the inverse map (\ref{inverse map}) there is
		\begin{align*}
		\hat{F}_{i\ell}^{*-1}\left(p+(-1)^{j+1}k_{n_i,1}\right)=\phi_p\left(\hat{F}_{i\ell}^{*}-(-1)^{j+1}k_{n_i,j}\right)
		\end{align*}
	From (\ref{bootmultidonsker}), one can infer with Slutzky's Lemma \citep[cf.][Example 1.4.7]{vaart_weak_2000} the following convergence in $D(\br)^{2d}$ for every $i\in\{1,\dots,k\}$
		\begin{align*}
		\sqrt{n_i}\left(\hat{F}_{i\ell}^*-k_{n_i,1}-\hat{F}_{i\ell},\hat{F}_{i\ell}^*+k_{n_i,2}-\hat{F}_{i\ell}\right)_{\ell\in\{1,\dots,d\}}\stackrel{d}{\longrightarrow}\left(\sqrt{\kappa_i}\mathbb{G}_i-\xi,\sqrt{\kappa_i}\mathbb{G}_i+\xi\right)
		\end{align*}
	As described in the proof of Proposition \ref{prop1}, the function $\phi_p$ is Hadamard differentiable and the functional delta method for the bootstrap \citep[cf.][Thm. 3.9.11]{vaart_weak_2000} can be applied with the function $\phi(G_1,G_2)=\left(\phi_p(G_1),\phi_p(G_2)\right)$ for every $\ell\in\{1,\dots,d\}$ and every $i\in\{1,\dots,k\}$:
		\begin{align*}
		&\sqrt{n_i}\left(\hat{F}_{i\ell}^{*-1}\left(p+(-1)^{j+1}k_{n_i,j}\right)-\hat{F}_{i\ell}^{-1}(p)\right)_{j\in\{1,2\}}\\
		=\quad&\sqrt{n_i}\left(\phi_p\left(\hat{F}_{i\ell}^*-(-1)^{j+1}k_{n_i,j}\right)-\phi_p\left(\hat{F}_{i\ell}\right)\right)_{j\in\{1,2\}}\\
		\stackrel{d}{\longrightarrow}&-\frac{1}{f_{i\ell}(q_{i\ell})}(\sqrt{\kappa_i}\mathbb{G}_i\left(q_{i\ell})+(-1)^j\xi\right)_{j\in\{1,2\}}.
		\end{align*}
	By applying this formula to the representation (\ref{formularbootint}) of the estimator, I get in probability
		\begin{align*}
		&\sqrt{n_i}\left(\hat{F}_{i\ell}^{*-1}\left(p+k_{n_i,1}\right)-\hat{F}_{i\ell}^{*-1}\left(p-k_{n_i,2}\right)\right)\\
		=\quad&\sqrt{n_i}\left(\hat{F}_{i\ell}^{*-1}\left(p+k_{n_i,1}\right)-\hat F_{i\ell}^ {-1}(p)\right)-\sqrt{n_i}\left(\hat{F}_{i\ell}^{*-1}\left(p-k_{n_i,2}\right)-\hat F_{i\ell}^{-1}(p)\right)\\
		\stackrel{P}{\longrightarrow}&\quad\cfrac{2z_{\alpha/2}\sqrt{p(1-p)}}{f_{i\ell}(q_{i\ell})}
		\end{align*}
	and, thus,
	\begin{align*}
	\hat\sigma_{i\ell}^{PB*}(p)\stackrel{P}{\longrightarrow}\frac{1}{2z_{\alpha/2}}\cfrac{2z_{\alpha/2}\sqrt{p(1-p)}}{f_{i\ell}(q_{i\ell})}=\cfrac{\sqrt{p(1-p)}}{f_{i\ell}(q_{i\ell})}.
	\end{align*}
	\end{proof}
From the two previous lemmas it follows, that $\hat{\mathbf{\Sigma}}^{PB*}\stackrel{P}{\longrightarrow}\mathbf{\Sigma}$ given the data.
The consistency of the bootstrap kernel density estimator is given analogously to Lemma S.3 of the supplement of \citet{ditzhaus_qanova_2021}.
	\begin{lemma}\label{bootkernel density}
	Let $\hat f_{K,i,\ell}^*$ be the bootstrap version of the kernel density estimator $\hat f_{K,i,\ell}$ in (\ref{kernel density}).
	Then, under Assumption \ref{kernelassumption} it holds the following  conditional convergence given the data in probability:
		\begin{align*}
		\sup_{x\in\br}\left|\hat f_{K,i,\ell}^*(x)-f_{i\ell}(x)\right|\stackrel{p}{\longrightarrow}0,\,i\in\{1,\dots,k\},\ell\in\{1,\dots,d\}.
		\end{align*}
	\end{lemma}
	\begin{proof}
	The proof runs analogously to the proof of Lemma S.3 in the supplement of \citet[p.~14]{ditzhaus_qanova_2021}.
	The lemma contains the assertion for an estimator based on permutation but this does not influence any argument of the proof.
	Instead of the pooled density $f$ in \citet{ditzhaus_qanova_2021}, the single densities $f_{i\ell},\,i\in\{1,\dots,k\},\ell\in\{1,\dots,d\}$ are considered.
	That is why further details are omitted.
	\end{proof}
From the previous lemma follows with Lemma \ref{bootjointconsistent} that $\hat{\mathbf{\Sigma}}^{K*}\stackrel{P}{\longrightarrow}\mathbf{\Sigma}$ given the data.
The consistency of the bootstrap version of the bootstrap plug-in estimator is not proven in \citet{ditzhaus_qanova_2021}.
Therefore, the proof will be omitted here as well.

\section{Proofs}\label{proofs}
The following section presents the proofs of the propositions and theorems of the paper.
	
	\subsubsection*{Proof of Proposition \ref{prop1}}
	Applying the delta-method for metrizable topological vector spaces \citep[cf.][Theorem 3.9.4]{vaart_weak_2000} to \eqref{multidonsker} yields the assertion.
	Let $\mathbb{D}=\{G:\br\to\br\,|\textrm{nondecreasing},\;\textrm{right- continuous}\}$, it holds $\mathbb{D}\subset D(\br)$ \citep[cf.][Supplement, p.~10]{ditzhaus_qanova_2021}.
	The function applied in the delta-method is the inverse mapping \citep[cf.][p.~385]{vaart_weak_2000}:
		\begin{align}\label{inverse map}
		\phi_p:\mathbb{D}\subset D(\br)\to\br,\;\phi_p(G)=G^{-1}(p)=\inf\{t\in\br|G(t)\ge p\},\;p\in(0,1).
		\end{align}
	For a fixed $\ell\in\{1,\dots,d\}$ the function $F_{i\ell}$ is in $\mathbb{D}$ and by Assumption \ref{ass} $F_{i\ell}$ is differentiable at $q_{i\ell}=F_{i\ell}^{-1}(p)$ with positive derivative $f_{i\ell}(q_{i\ell})$.
	Thus, by Lemma 3.9.20 of \citet[p.~385]{vaart_weak_2000} the function $\phi_p$ is for every $p\in(0,1)$ Hadamard differentiable at $F_{i\ell}$ tangentially to the space $\mathbb{D}_{q_{i\ell}}\subset D(\br)$, which contains all function $\alpha\in\mathbb{D}(\br)$ that are continuous at $q_{i\ell}$.
	The Hadamard derivative is calculated by
		\[\phi_{p,F_{i\ell}}'(\alpha)=-\frac{\alpha(q_{i\ell})}{f_{i\ell}(q_{i\ell})}.\]
	It holds $\sqrt{n_i}\xlongrightarrow{n_i\to\infty}\infty$ and $\mathbb{G}_i$ is separable.
	Otherwise the separable version of $\mathbb{G}_i$ is chosen as described in \citet[Section 2.2.3]{vaart_weak_2000}.
	Thus, the requirements of Theorem 3.9.4 in \citet{vaart_weak_2000} are fulfilled and it follows that
		\begin{align*}
		\sqrt{n_i}\left(\phi_p(\hat{F}_{i\ell})-\phi_p(F_{i\ell})\right)&=\sqrt{n_i}\left(\hat{F}_{i\ell}^{-1}(p)-F_{i\ell}^{-1}(p)\right)\\
		&=\sqrt{n_i}\left(\hat{q}_{i\ell}-q_{i\ell}\right)
		\end{align*}
	converges for every $\ell\in\{1,\dots,d\}$ weakly to
		\begin{align*}
		\phi'_{F_{i\ell}}(\mathbb{G}_i)&=-\frac{\mathbb{G}_i\left(q_{i\ell}\right)}{f_{i\ell}(q_{i\ell})}.
		\end{align*}
	For $a,b\in\{1,\dots,d\}$ it holds
		\[E\left(-\frac{\mathbb{G}_i\left(q_{ia}\right)}{f_{ia}(q_{ia})}\right)=-\frac{1}{f_{ia}(q_{ia})}E\left(\mathbb{G}_i\left(q_{ia}\right)\right)=0\]
	and
		\begin{align}
		&Cov\left(-\frac{\mathbb{G}_i\left(q_{ia}\right)}{f_{ia}(q_{ia})},-\frac{\mathbb{G}_i\left(q_{ib}\right)}{f_{ib}(q_{ib})}\right)\nonumber\\
		=&\frac{1}{f_{ia}(q_{ia})f_{ib}(q_{ib})}Cov\left(\mathbb{G}_i\left(q_{ia}\right),\mathbb{G}_i\left(q_{ib}\right)\right)\nonumber\\
		=&\begin{cases}
		\cfrac{1}{f_{ia}^2(q_{ia})}\left[F_{ia}\left(\min(q_{ia},q_{ia})\right)-F_{ia}(q_{ia})F_{ia}(q_{ia})\right],\;a=b\\
		\cfrac{1}{f_{ia}(q_{ia})f_{ib}(q_{ib})}\left[F_{iab}(q_{ia},q_{ib})-F_{ia}(q_{ia})F_{ib}(q_{ib})\right],\;a\not=b
		\end{cases}\nonumber\\
		=&\begin{cases}
		\cfrac{1}{f_{ia}^2(q_{ia})}\left[F_{ia}\left(F_{ia}^{-1}(p)\right)-F_{ia}\left(F_{ia}^{-1}(p)\right)F_{ia}\left(F_{ia}^{-1}(p)\right)\right],\;a=b\\
		\cfrac{1}{f_{ia}(q_{ia})f_{ib}(q_{ib})}\left[F_{iab}\left(F_{ia}^{-1}(p),F_{ib}^{-1}(p)\right)-F_{ia}\left(F_{ia}^{-1}(p)\right)F_{ib}\left(F_{ib}^{-1}(p)\right)\right],\;a\not=b
		\end{cases}\nonumber\\
		=&\begin{cases}
		\cfrac{1}{f_{ia}^2(q_{ia})}\left[p-p^2\right],\;a=b\\
		\cfrac{1}{f_{ia}(q_{ia})f_{ib}(q_{ib})}\left[F_{iab}\left(F_{ia}^{-1}(p),F_{ib}^{-1}(p)\right)-p^2\right],\;a\not=b
		\end{cases}.
		\end{align}
	Consequently, this means $\frac{\sqrt{n_i}}{\sqrt{n}}\,\sqrt{n}\left(\hat{q}_{i\ell}-q_{i\ell}\right)=\sqrt{n_i}\left(\hat{q}_{i\ell}-q_{i\ell}\right)\stackrel{d}{\to}\sqrt{\kappa_i}\,\mathbf{Z}_i$ and with the assumption of non-vanishing groups the assertion follows.
		
	\subsubsection*{Proof of Theorem \ref{distributiontests}}
	\begin{enumerate}
	\item
	The statement from Proposition \ref{prop1} is for every $i\in\{1,\dots,k\}$:
		\[\sqrt{n}\left(\mathbf{\hat q}_i-\mathbf{q}_i\right)\stackrel{d}{\longrightarrow}\mathbf{Z}_i\sim\mathcal{N}\left(\mathbf{0}_{d},\mathbf{\Sigma}_i\right)\]
	 and from the independent distributed groups $i\in\{1,\dots,k\}$ it follows
		\[\sqrt{n}\left(\mathbf{\hat q}-\mathbf{q}\right)\stackrel{d}{\longrightarrow}\mathbf{Z}\sim\mathcal{N}\left(\mathbf{0}_{dk},\mathbf{\Sigma}\right).\]
	Under $\mathcal{H}_0$ and by the Continuous Mapping Theorem \citep[Thm. 1.11.1]{vaart_weak_2000} the following is also true:
		\[\sqrt{n}\mathbf{T\hat q}=\sqrt{n}\left(\mathbf{T\hat q}-\mathbf{Tq}\right)=\mathbf{T}\sqrt{n}\left(\mathbf{\hat q}-\mathbf{q}\right)\stackrel{d}{\longrightarrow}\mathbf{TZ}\sim\mathcal{N}\left(\mathbf{0}_{dk},\mathbf{T\Sigma T}\right).\]
	Therefore, due to the consistency of $E(\mathbf{T},\mathbf{\hat\Sigma})$ and by Slutzky's theorem \citep[cf.][Example 1.4.7]{vaart_weak_2000} it follows
		\begin{align*}
		S_n(\mathbf{T})&=n(\mathbf{T\hat q})'E\left(\mathbf{T},\mathbf{\hat\Sigma}\right)\mathbf{T\hat q}\\
		&=(\sqrt{n}\mathbf{T\hat q})'E\left(\mathbf{T},\mathbf{\hat\Sigma}\right)\left(\sqrt{n}\mathbf{T\hat q}\right)\stackrel{d}{\to}(\mathbf{TZ})'E\left(\mathbf{T},\mathbf{\Sigma}\right)\left(\mathbf{TZ}\right).
		\end{align*}
	By Theorem 8.35 in \citet{brunner_rank-_2018} this has the same distribution as the random variable $\sum_{i=1}^{dk}\lambda_iB_i$ with $B_i\stackrel{iid}{\sim}\chi_1^2,\,i\in\{1,\dots,dk\}$, and $\lambda_i$ are the eigenvalues of $E\left(\mathbf{T},\mathbf{\Sigma}\right)\mathbf{T\Sigma T}$.
	\item
	Proposition \ref{prop1} is not restricted to the null hypothesis. 
	Thus, it can be concluded from the Continuous Mapping Theorem \citep[Thm. 1.11.1]{vaart_weak_2000} that $n^ {-1}S_n(\mathbf{T})$ always converges in probability to $(\mathbf{Tq})'E\left(\mathbf{T},\mathbf{\Sigma}\right)\mathbf{Tq}$.
	That is why it remains to prove that from $\mathcal{H}_1:\,\mathbf{Tq}\not=\mathbf{0}_{dk}$ it always follows that $(\mathbf{Tq})'E\left(\mathbf{T},\mathbf{\Sigma}\right)\mathbf{Tq}>0$.
	Let $\mathbf{Tq}\not=\mathbf{0}_{dk}$.
	I need to consider the two versions of $E\left(\mathbf{T},\mathbf{\Sigma}\right)$ separately.
	The proof of the ATS follows immediately with
		\begin{align*}
		&\mathbf{Tq}\not=\mathbf{0}_{dk}\Rightarrow(\mathbf{Tq})'\mathbf{Tq}>0\Rightarrow\frac{(\mathbf{Tq})'\mathbf{Tq}}{\tr(\mathbf{T\Sigma T})}>0
		\end{align*}
	The proof of the MATS follows analogously  to the proof of Theorem 2 in \citet{ditzhaus_qanova_2021}.
	The covariance matrix $\mathbf{\Sigma}$ is positive semidefinite and symmetric by definition.
	Thus, the square root $\mathbf{\Sigma}^{\frac{1}{2}}$ exists and is also positive semidefinite and symmetric \citep[Thm. 7.2.6]{horn_matrix_2013}.
	As a consequence, the square root $\mathbf{\Sigma}_0^{\frac{1}{2}}$ exists as well.
	Moreover, there is some $\tilde{\mathbf{q}}\in\br^{dk}$ such that $\mathbf{q}=\mathbf{\Sigma}_0^{\frac{1}{2}}\tilde{\mathbf{q}}$.
	From the preceding and the non-singularity of $\mathbf{\Sigma}_0$ it follows
		\begin{align*}
		\mathbf{T\Sigma}_0^{\frac{1}{2}}\left[\left(\mathbf{T\Sigma}_0^{\frac{1}{2}}\right)^+\mathbf{Tq}\right]=&\;\mathbf{T\Sigma}_0^{\frac{1}{2}}\left[\left(\mathbf{T\Sigma}_0^{\frac{1}{2}}\right)^+\mathbf{T\Sigma}_0^{\frac{1}{2}}\tilde{\mathbf{q}}\right]\\
		=&\;\mathbf{T\Sigma}_0^{\frac{1}{2}}\tilde{\mathbf{q}}=\mathbf{T\Sigma}_0^{\frac{1}{2}}\left(\mathbf{\Sigma}_0^{\frac{1}{2}}\right)^+\mathbf{q}\\
		=&\;\mathbf{Tq}\not=\mathbf{0}_{dk}.
		\end{align*}
	And with the well known properties $\left(\mathbf{A}'\right)^+=\left(\mathbf{A}^+\right)'$ and $\left(\mathbf{A}'\mathbf{A}\right)^+=\mathbf{A}^+\left(\mathbf{A}'\right)^+$ of Moore-Penrose inverses \citep[cf.][p.~67]{rao_generalized_1971} I conclude
		\begin{align*}
		\left(\mathbf{Tq}\right)'\left(\mathbf{T}\mathbf{\Sigma}_0\mathbf{T}\right)^+\mathbf{Tq}=&\left(\mathbf{Tq}\right)'\left[\left(\mathbf{\Sigma}_0^{\frac{1}{2}}\mathbf{T}\right)'\mathbf{\Sigma}_0^{\frac{1}{2}}\mathbf{T}\right]^+\mathbf{Tq}\\
		=&\left(\mathbf{Tq}\right)'\left(\mathbf{\Sigma}_0^{\frac{1}{2}}\mathbf{T}\right)^+\left(\mathbf{T}\mathbf{\Sigma}_0^{\frac{1}{2}}\right)^+\mathbf{Tq}\\
		=&\left(\mathbf{Tq}\right)'\left[\left(\mathbf{T}\mathbf{\Sigma}_0^{\frac{1}{2}}\right)^+\right]'\left(\mathbf{T}\mathbf{\Sigma}_0^{\frac{1}{2}}\right)^+\mathbf{Tq}\\
		=&\left[\left(\mathbf{T}\mathbf{\Sigma}_0^{\frac{1}{2}}\right)^+\mathbf{Tq}\right]'\left[\left(\mathbf{T}\mathbf{\Sigma}_0^{\frac{1}{2}}\right)^+\mathbf{Tq}\right]>0.
		\end{align*}
	\end{enumerate}		
	\subsubsection*{Proof of Proposition \ref{jointconsistent}}
	Even if \citet{babu_joint_1988} proved this, the result can be easily shown with the theory of empirical processes.
	By Example 2.1.3 in \citet{vaart_weak_2000} the $2$-dimensional empirical distribution functions $F_{i\ell m}$ can be identified with the empirical measure indexed by $\mathcal{F}=\{\mathbf{1}_{(-\infty,t_\ell]\times(-\infty,t_m]}|(t_\ell,t_m)\in\br^2\}$, which forms a Donsker-Class for $i\in\{1,\dots,k\},\,\ell,m\in\{1,\dots,d\}$.
	Thus, the set $\mathcal{F}$ is also a Glivenko-Cantelli-Class \citep[cf.][Lemma 2.4.5]{vaart_weak_2000}:
		\[\sup_{(t_1,t_2)\in\br^2}\left|\hat F_{i\ell m}(t_1,t_2)- F_{i\ell m}(t_1,t_2)\right|\stackrel{a.e.}{\longrightarrow}0.\]
	From the continuity of $F_{im\ell}$ at $(q_{i\ell},q_{im})$ and $(\hat q_{i\ell},\hat q_{im})\stackrel{a.e.}{\longrightarrow}(q_{i\ell},q_{im})$ it follows from this analogous to \citet[cf.][p.~18]{babu_joint_1988}:
		\[\hat F_{i\ell m}\left(\hat q_{i\ell},\hat q_{im}\right)\stackrel{a.e.}{\longrightarrow}F_{i\ell m}(q_{i\ell},q_{im}).\]
	This yields the assertion.
\subsubsection*{Proof of Theorem \ref{bootdistributiontest}}
The proof is analogous to the proof of Theorem \ref{distributiontests}.
The bootstrap version of \eqref{multidonsker} follows from Theorem 3.6.2, Chapter 1.12 in \citet{vaart_weak_2000} and from Lemma \ref{VC}:
	\begin{align}\label{bootmultidonsker}
	\sqrt{n_i}\left(\hat{F}_{i\ell}^*-\hat{F}_{i\ell}\right)_{\ell\in\{1,\dots,d\}}\stackrel{d}{\longrightarrow}\mathbb{G}_i\;\text{in}\;D(\br)^d,\;i\in\{1,\dots,k\},
	\end{align}
given the data in probability.
Here, $\hat{F}_{i\ell}^*$ describes the empirical distribution function calculated with the bootstrap sample, e.g. $\hat F_{i\ell}^*(t)=\frac{1}{n_i}\sum_{j=1}^{n_i}\mathbf{1}_{\{X_{ij\ell}^*\le t\}}$.
The delta-method for bootstrapping \citep[Theorem 3.9.11]{vaart_weak_2000} yields analogous to Proposition \ref{prop1} under Assumption \ref{non-vanishing}:
	\begin{align}\label{bootprop1}
 	\sqrt{n}\left(\hat q_{i\ell}^*-\hat q_{i\ell}\right)_{\ell\in\{1,\dots,d\}}\stackrel{d}{\longrightarrow}\mathbf{Z}_i,\;i\in\{1,\dots,k\},
 	\end{align}
given the data in probability, where $\hat q_{i\ell}^*=F_{i\ell}^{*-1}(p)$ is the empirical bootstrap quantile and $\mathbf{Z}_i$ is as in Proposition \ref{prop1}.
As a result, the bootstrap version of the central limit theorem gives us the same limit process $\mathbf{Z}_i,\,i\in\{1,\dots,k\}$ as the regular one.
That is why the covariance of the limit process is again $\mathbf{\Sigma}=\bigoplus_{i=1}^k\mathbf{\Sigma}^{(i)}$ and it is needed to estimate $\mathbf{\Sigma}$ with the bootstrap sample.
	By \eqref{bootprop1} it holds:
		\[\sqrt{n}\left(\mathbf{\hat q}^*-\mathbf{\hat q}\right)\stackrel{d}{\longrightarrow}\mathbf{Z}\sim\mathcal{N}\left(\mathbf{0}_{dk},\mathbf{\Sigma}\right).\]
	From the Continuous Mapping Theorem \citep[Thm. 1.11.1]{vaart_weak_2000} it follows
		\[\sqrt{n}\mathbf{T}\left(\mathbf{\hat q}^*-\mathbf{\hat q}\right)\stackrel{d}{\longrightarrow}\mathbf{TZ}\sim\mathcal{N}\left(\mathbf{0}_{dk},\mathbf{T\Sigma T}\right).\]
	The consistency of $E(\mathbf{T},\mathbf{\hat\Sigma}^*)$ follows from Proposition \ref{consistent}.
	Identical to the proof of the Theorem \ref{distributiontests}, this yields the assertion.

\end{document}